\documentclass[a4paper,11pt,bibliography=totoc,listof=totoc,headinclude=true,cleardoublepage=empty,oneside]{article}

\usepackage[english]{babel}
\usepackage{tikz}
\usepackage[utf8]{inputenc}
\usepackage{ifthen}
\usepackage{color}
 \usepackage{standalone}
\usepackage{amsmath,amsthm,amssymb,amsfonts}
\usepackage{graphicx}
\usepackage{psfrag}
\usepackage{pgfplots}

\usepackage[margin=2.5 cm]{geometry}
\usetikzlibrary{calc}

\usepackage[unicode,colorlinks=true,linkcolor=blue,citecolor=green,urlcolor=black,pagebackref=false]{hyperref}

\definecolor{change}{rgb}{0,.55,.55}

\newcommand{\R}{\mathbb{R}}
\newcommand{\N}{\mathbb{N}}
\newcommand{\SN}{\mathcal{N}}
\renewcommand{\S}{\mathcal{S}}

\newcommand{\T}{\mathcal{T}}
\newcommand{\diff}{\mathrm{d}}

\newcommand{\TR}{T_{\text{ref}}}

\newcommand{\omegaref}{\omega_{\text{ref}}}
\newcommand{\Pat}{\operatorname{Pat}(T)}

\newcommand{\SCP}[3]{\langle #1, #2 \rangle _{#3}}
\newcommand{\norm}[2]{\left\|#1\right\|_{#2}}
\newcommand{\Tunif}{\T_{\ell,u}}

\newcommand{\vertiii}[1]{{\left\vert\kern-0.25ex\left\vert\kern-0.25ex\left\vert #1
    \right\vert\kern-0.25ex\right\vert\kern-0.25ex\right\vert}}
\newenvironment{fatproof}{{\bfseries Proof.}}{\qed \newline}

\def\restrict#1{\raise-.5ex\hbox{\ensuremath|}_{#1}}

\newcommand{\exspace}[1]{\widetilde{H}^s(#1)}

\newtheorem{theorem}{Theorem}[section]
\newtheorem{lemma}[theorem]{Lemma}
\newtheorem{algo}[theorem]{Algorithm}
\theoremstyle{definition}
\theoremstyle{definition}\newtheorem{remark}[theorem]{Remark}

\usetikzlibrary{shapes.geometric, arrows}
\tikzstyle{boxs} = [rectangle, rounded corners, minimum width=2cm, minimum height=1cm,text centered, draw=black]
\tikzstyle{arrow} = [thick,->,>=stealth]

\usepackage[section]{placeins}

\numberwithin{equation}{section}

\title{Two-level error estimation for the integral fractional Laplacian}

\author{Markus Faustmann, Ernst Peter Stephan and David Wörgötter}
\date{\today}

\begin{document}
\maketitle
\begin{abstract}
For the singular integral definition of the fractional Laplacian, we consider an adaptive finite element method steered by two-level error indicators. For this algorithm, we show linear convergence in two and three space dimensions as well as convergence of the algorithm with optimal algebraic rates in 2D, when newest vertex bisection is empolyed for mesh refinement.
\end{abstract}

           \section{Introduction}
Solutions to PDEs with non-integer powers of differential operators, such as the fractional Laplacian $(-\Delta)^s$, $s \in (0,1)$, which are commonly used to describe anomalous, non-local diffusion processes, typically exhibit a singular behaviour at the boundary of the computational domain.
Thus, solutions to fractional PDEs are, in general, less regular than in the case of integer order powers, \cite{Grubb15} and as a consequence finite element approximations on uniform meshes lead to non-optimal orders of convergence, \cite{AB17}. 

The nature of this boundary layer singularity is well understood, \cite{Grubb15,GSS21}, and regularity estimates with weights can be obtained, \cite{AB17,BN21,FMMS21}, which in turn can be used to a-priori design suitable meshes that allow to recuperate better convergence rates. In this sense, graded meshes with appropriate grading factors give algebraic rates for the finite regularity case, \cite{AB17,GSS21}, while exponentially graded meshes can be used to obtain exponentially convergent $hp$-FEM approximations in the case of weighted analytic regularity, \cite{FMMS22}.

A different approach to dealing with the singularities at the boundary, which we employ here, is to use adaptive finite element methods (AFEM). Hereby, meshes are locally refined only where some error measure indicates large errors. We consider a classical adaptive loop of the form  {\texttt{SOLVE}--\texttt{ESTIMATE}--\texttt{MARK}--\texttt{REFINE}}, see, e.g, \cite{Dorfler}. The key to the success of such an adaptive strategy is the choice of error indicators for the \texttt{ESTIMATE} step. 

A classical strategy for a-posteriori error estimation is to take the (weighted) residual, which in our case takes the form $\|f-(-\Delta)^s u_h\|_{L^2}$, with $f$ denoting the data and $u_h$ being the FEM solution, as local error indicator. After correcting this residual with a suitable weight function, such that one indeed obtains a well-defined $L^2$-function, one obtains an algorithm that converges with optimal algebraic rates, \cite{FMP21a}. For other strategies on a-posteriori error estimation for classical, integer order differential operators, we refer to \cite{BW85,ZZ87,MS99,CFPP14}.
While the weighted residual error estimator has nice analytical properties, it has distinct drawbacks for fractional PDEs. Most notably, the computation of the $L^2$-norm enforces the use of accurate, expensive quadrature rules as the integrand tends to be singular, and the evaluation of the fractional Laplacian in the quadrature points is expensive as well. 

Therefore, in the following, we consider a different strategy for error estimation using approximations on different refinement levels to measure the local errors. The advantage of this approach is that the costly evaluation of the fractional differential operator is avoided. Two level or hierarchical error estimators are commonly used in practice, \cite{MMS97,MSW98,EFFS09,stephan2} and perform well in numerical examples. On the analytical side, reliability of the error indicators always hinges on certain saturation assumptions, i.e, that refinement reduces the error by a contraction factor $q<1$. Under these assumptions, there even holds optimality for adaptive finite element and boundary element methods, \cite{PRS20}.
In this article, we show a corresponding result for the fractional Laplacian. 
In fact, for $d=2,3$ and under certain saturation assumptions, we show discrete efficiency and reliability of the error indicators, which implies linear convergence of the error between exact solution and AFEM approximation. For the case of two spatial dimensions, we show discrete stability as well, which together with the frameworks of \cite{CFPP14,PRS20} even provides optimal algebraic convergence of the algorithm.   

The present paper is structured as follows: In Section~2, we introduce the model problem, its discretization as well as the adaptive finite element method based on two-level error estimation. Finally, we present our main results, Theorem~\ref{thm:mainresult1}-\ref{thm:mainresult3}. Section~3 is dedicated to the proofs of the main results. A key step that might be of independent interest is the equivalence of nodal interpolation and Scott-Zhang projection on discrete spaces in the $L^2$-sense, Theorems~\ref{thm:NodalSZ}. Finally, Section~4 presents some numerical examples that underline the optimal convergence of the adaptive algorithm.

Throughout this article, we write $\lesssim$ to abbreviate $\leq$ up to a generic constant $C>0$ that does not depend on critical parameters in our analysis. Moreover, we write $\simeq$ to indicate that both estimates $\lesssim$ and $\gtrsim$ hold.

	\section{Setting and main results}

Let $\Omega\subseteq \R^d$ with $d\in \{2,3\}$ be a bounded polygonal/polyhedral Lipschitz domain and $s\in(0,1)$. We note that there are several ways to define non-integer powers of differential operators such as the fractional Laplacian $(-\Delta)^s$. A classical definition on the full space is given as operator with Fourier symbol $|\zeta|^{2s}$, but definitions using spectral theory, semi-group approaches or PDE-extensions can also be employed, \cite{CS07,fraclapintro,tendef}. We note that these approaches are equivalent on the full space, but not on a bounded domain. 

Here, for a sufficiently smooth function $u$ defined on $\R^d$ and $x\in\Omega$, we use the integral fractional Laplacian $(-\Delta)^s u$, defined pointwise as the Cauchy principal value
\begin{align*}
	(-\Delta)^su(x) := C(d,s)\ \mathrm{P.V.}\int_{\R^d}\frac{u(x)-u(y)}{|x-y|^{d+2s}}\ \diff y, \hskip 0.5 cm C(d,s) := \frac{2^{2s}s\ \Gamma\left(\frac{d}{2}+s\right)}{\pi^{d/2}\Gamma(1-s)},
\end{align*}
where $\Gamma(\cdot)$ is the Gamma function.

In this paper, we are interested in the solution $u$ of the Dirichlet problem
\begin{align*}
	(-\Delta)^s u &= f \hskip 0.5cm \mathrm{in}\ \Omega\\ 
	u &= 0 \hskip 0.53cm \mathrm{in}\ \R^d\setminus\overline{\Omega} ,
\end{align*}
where $f\in L^ 2(\Omega)$ is a given right-hand side. The proper function spaces for this problem are fractional Sobolev spaces with an exterior boundary condition defined by means of the 
Aronstein-Slobodeckij norm for any open set $\omega \subset \R^d$
\begin{align*}
|v|^2_{H^s(\omega)} &= \int_{x \in \omega} \int_{y \in \omega} \frac{|v(x) - v(y)|^2}{|x-y|^{d+2s}}
\,\diff x\,\diff y, 
\qquad \|v\|^2_{H^s(\omega)} = \|v\|^2_{L^2(\omega)} + |v|^2_{H^s(\omega)},  \\
 \widetilde{H}^{s}(\Omega) &:= \{u \in H^s(\R^d)  \,| \, u\equiv 0 \; \text{on} \; \R^d \backslash \overline{\Omega} \}, \quad \|v\|^2_{\widetilde H^s(\Omega)} = \|v\|^2_{H^s(\Omega)} + \|v/r_{\partial\Omega}^s\|^2_{H^s(\Omega)},
\end{align*}
where $r_{\partial\Omega}(x) := \operatorname*{dist}(x,\partial\Omega)$ is the Euclidean distance of a point $x \in \Omega$ to the boundary of $\Omega$.

The weak formulation of our model problem reads as: Find $u\in \exspace{\Omega}$ such that
\begin{align}\label{setting:weakform}
a(u,w) :=\frac{C(d,s)}{2}\int_{\R^d}\int_{\R^d}\frac{(u(x)-u(y))(w(x)-w(y))}{|x-y|^{d+2s}} \ \mathrm{d}y\mathrm{d}x = \int_{\Omega}fw\ \diff x
\end{align}
for all $w\in \exspace{\Omega}$, \cite{AB17}. For $u\in\exspace{\Omega}$, the energy-norm $\vertiii{u} := a(u,u)^{1/2}$ is an equivalent norm on $\exspace{\Omega}$ (see \cite{AB17}), thus, for given $f\in L^2(\Omega)$, the weak formulation \eqref{setting:weakform} has a unique solution.

\subsection{Discretization}
In the following, we consider a finite element discretization of the fractional PDE, which is 
based on regular (in the sense of Ciarlet, \cite{Ciarlet78}, i.e., there are no hanging nodes) meshes $\T_{\ell}$ that decompose $\Omega$ into (open) triangles/tetrahedrons. The subscript $\ell\in\N_0$ will refer to the $\ell$-th step of an iterative algorithm. 
Moreover, we assume that all triangulations employed are $\gamma$-shape regular, meaning that there exists a constant $\gamma > 0$, such that
$$
\max_{T\in\T_{\ell}}\left(\operatorname{diam}(T)^d/|T|^{1/d}\right)\leq \gamma,
$$
where $\operatorname{diam}(T)$ denotes the Euclidean diameter of a triangle $T$ and $|T|$ denotes the Lebesgue-measure of $T$. 

\medskip 

The set of all nodes of $\T_{\ell}$ (i.e., the set of all vertices of elements in the triangulation $\T_{\ell}$) that lie inside of $\Omega$ and not on the boundary $\partial \Omega$ is denoted by $\mathcal{N}_{\T_{\ell}}^*$. Furthermore, for $T\in\T_{\ell}$, we write $\mathcal{N}_{\T_{\ell}}^*(T)$ for the set of all vertices of $T$ that lie inside of $\Omega$. For an element $T\in\T_{\ell}$, the set 
$$
\Omega_{\ell}[T] := \operatorname*{interior}\bigcup\{\overline{T'} \ |\ T'\in\T_{\ell}, \overline{T}\cap \overline{T'}\neq \emptyset\}
$$ 
is called the element patch of $T$. 

For a triangulation $\T_{\ell}$, the associated mesh-width function $h_{\ell}\in L^{\infty}(\Omega)$ is elementwise defined by $h_{\ell}\restrict{T}:= |T|^{1/d}$ for $T\in\T_{\ell}$. 

\bigskip 

 Based on a regular triangulation $\T_{\ell}$, we define the space of globally continuous, piecewise linear polynomials as 
\begin{align*}
	\S^1(\T_{\ell}) &:= \{u\in C(\overline\Omega)\ |\ u\restrict{T}\in \mathcal{P}^1(T)\text{ for all }T\in\T_{\ell}\}, \\
	\S_0^1(\T_{\ell})&:=\S^1(\T_{\ell})\cap H_0^1(\Omega),
\end{align*}
where $\mathcal{P}^ 1(T)$ is the space of all linear polynomials on $T$. 

\bigskip

Due to $H^1(\Omega)\subseteq H^s(\Omega)$ for $s\in (0,1)$ (see \cite[Propositions 2.1 and 2.2]{frac_sob}), we obtain that $\S_0^1(\T_{\ell})$ is a closed subspace of $\exspace{\Omega}$ for all regular triangulations $\T_{\ell}$. The discretization of the weak formulation \eqref{setting:weakform} thus reads as: Find the Galerkin solution $u_{\ell}\in \S_0^1(\T_{\ell})$, which satisfies
\begin{align}\label{setting:discrweakform}
a(u_{\ell}, w_{\ell}) = \int_{\Omega}fw_{\ell}\ \diff x \qquad \forall w_{\ell}\in \S_0^1(\T_{\ell}).
\end{align}
By the Lax-Milgram lemma, one also obtains unique solvability of the discrete Galerkin formulation.
Henceforth, by $u_{\ell}$ we always denote the piecewise linear Galerkin solution corresponding to $\T_{\ell}$, whereas $u$ shall always refer to the exact solution of \eqref{setting:weakform}.

\subsection{Mesh refinement}\label{pattern}

In this work we consider meshes obtained from consecutive refinement of an initial triangulation $\T_0$. We assume that refinement is done by newest vertex bisection (NVB). That is, bisecting an element $T$ of a triangulation $\T_{\ell}$ always follows the same principle, namely that an edge opposite to the newest vertex is halved, see e.g. \cite{NVB} for $d=2$ or \cite{Stevenson2} for $d\geq 2$.

Consider a triangulation $\T_{\ell}$ and a set $\mathcal{M}_{\ell}\subseteq \T_{\ell}$ of marked elements. As in \cite{PRS20, Bespalov}, we use the following notation: When $d=2$, the triangulation $\T_{\ell+1}:=\operatorname{refine}(\T_{\ell},\mathcal{M}_{\ell})$ shall denote the coarsest NVB refinement of $\T_{\ell}$ such that all edges of marked elements $T\in\T_{\ell}$ are bisected. This corresponds to bisecting every marked element three times into four sons, see e.g. \cite{Erath}.

\medskip 

Newest vertex bisection in 3D is a bit more complex: In order to guarantee $\gamma$-shape regularity of refined meshes, each element $T\in\T_{\ell}$ has a certain refinement type $\beta\in\{0,1,2\}$ which specifies the next bisection step, see \cite{Stevenson2, Erath, stephan1} for further discussion.

In the three-dimensional case, it was required in \cite{Erath} that $\T_{\ell+1}:=\operatorname{refine}(\T_{\ell},\mathcal{M}_{\ell})$ is the coarsest NVB refinement of $\T_{\ell}$ such that all faces of marked elements contain an interior node. According to \cite{Erath}, this can be done by splitting every marked element according to a certain rule, where the exact rule depends on the refinement parameter $\beta\in\{0,1,2\}$ of $T$. If $\beta=0$ or $\beta=1$, the element $T$ is split into $18$ sons with $14$ nodes. If $\beta=2$, then $T$ is split into $20$ sons with $14$ nodes. The splitting of a marked element $T$ is such that the $14$ new nodes always consist of the $4$ original nodes of $T$, the midpoints of the $6$ edges of $T$ and one node in the interior of each of the $4$ faces of $T$, respectively. See \cite{Erath} for a thorough discussion.

\medskip 

In this work we go a step further and require that, in 3D, the refinement $\T_{\ell+1}:=\operatorname{refine}(\T_{\ell},\mathcal{M}_{\ell})$ is the coarsest NVB refinement of $\T_{\ell}$ such that all faces of elements $T\in\Lambda(\mathcal{M}_{\ell})$ contain an interior node. Here, $\Lambda(\mathcal{M}_{\ell})$ is defined as
\begin{align*}
	\Lambda(\mathcal{M}_{\ell}) := \bigcup\left\{T'\in\T_{\ell}\ |\ T'\in\mathcal{M}_{\ell}\ \text{or\ } T'\  \text{shares\ at\ least\ one\ edge\ with\ an\  element\ }T\in\mathcal{M}_{\ell}\right\},
\end{align*}
i.e., it is the set of all elements in $\T_{\ell}$ that are either marked themselves or share one or more edges with a marked element. 
If we denote the refinement strategy from \cite{Erath} as $\operatorname{refine}_{\operatorname{EGP}}(\cdot)$ and the one used in this work as $\operatorname{refine}_{\operatorname{FSW}}(\cdot)$, then we observe that there holds 
\begin{align}\label{newrefmarked}
	\operatorname{refine}_{\operatorname{FSW}}(\T_{\ell}, \mathcal{M}_{\ell}) = \operatorname{refine}_{\operatorname{EGP}}(\T_{\ell}, \Lambda(\mathcal{M}_{\ell}))
\end{align}
for all subsets $\mathcal{M}_{\ell}\subseteq\T_{\ell}$.

\begin{remark}\label{setting:samesons}
	With this definition of a refinement rule, the refinement of a marked element does not depend on its neighbors, i.e., the mesh closure step does not further refine sons of marked elements. In particular, for refinement $\T_{\ell+1}:=\operatorname{refine}(\T_{\ell},\mathcal{M}_{\ell})$ and the uniform refinement $\Tunif:=\operatorname{refine}(\T_{\ell},\T_{\ell})$, there holds 
	\begin{align*}
		\{T\in\T_{\ell+1}\ |\ T\subseteq \mathcal{M}_{\ell}\} = \{T\in\Tunif\ |\ T\subseteq \mathcal{M}_{\ell}\},
	\end{align*}
	that is, the sons of marked elements are the same in $\T_{\ell+1}$ and $\Tunif$. 
	In the two-dimensional case this fact is easy to check, for $\operatorname{refine}_{\operatorname{EGP}}(\cdot)$ we refer to \cite{Erath}, and due to \eqref{newrefmarked}, the arguments provided in \cite{Erath} extend to $\operatorname{refine}_{\operatorname{FSW}}(\cdot)$. Moreover, if $\operatorname{refine}_{\operatorname{FSW}}(\cdot)$ is used, the sons of all elements $T\in\Lambda(\mathcal{M}_{\ell})$ are the same in $\T_{\ell+1}$ and $\Tunif$.
\end{remark}

\begin{remark}
	Note that for a marked element $T\in\mathcal{M}_{\ell}$, the number of elements $T'\in\T_{\ell}$ that share an edge with $T$ is bounded in terms of the $\gamma$-shape regularity of $\T_{\ell}$. That is, there holds 
	\begin{align*}
	\#\mathcal{M}_{\ell}\leq \#\Lambda(\mathcal{M}_{\ell})\leq C\ \#\mathcal{M}_{\ell},
	\end{align*}
	with $C>0$ depending only on the $\gamma$-shape regularity of $\T_{\ell}$. That is, the number of elements in $\Lambda(\mathcal{M}_{\ell})$ can be controlled in terms of the number of elements in $\mathcal{M}_{\ell}$ and this suggests that it might be possible to control the amount of refined elements in $\operatorname{refine}_{\operatorname{FSW}}(\T_{\ell},\mathcal{M}_{\ell})$ in terms of the number of refined elements in $\operatorname{refine}_{\operatorname{EGP}}(\T_{\ell},\mathcal{M}_{\ell})$.
\end{remark}

\medskip 

From now on, if not stated otherwise, we will only consider $\operatorname{refine}_{\operatorname{FSW}}(\cdot)$ as refinement strategy and we write $\operatorname{refine}(\cdot)$ instead of $\operatorname{refine}_{\operatorname{FSW}}(\cdot)$.

\medskip

Let us mention that, for a one-level refinement $\T_{\ell+1}:=\operatorname{refine}(\T_{\ell},\mathcal{M}_{\ell})$ of $\T_{\ell}$, elementary geometric considerations provide equivalence of the corresponding mesh-size functions 
\begin{align}\label{setting:geometric}
	h_{\ell+1}\leq h_{\ell}\leq C \ h_{\ell+1},
\end{align}
where $C=2$ for $d=2$ and $C=\sqrt[3]{32}$ for $d=3$.

\bigskip

Remark \ref{setting:samesons} implies the following lemma:
\begin{lemma}\label{setting:uniforminref}
	 For a triangulation $\T_{\ell}$ and a subset $\mathcal{M}_{\ell}\subseteq \T_{\ell}$ consider the refinement $\T_{\ell+1}:=\operatorname{refine}(\T_{\ell},\mathcal{M}_{\ell})$, as well as the uniform refinement $\Tunif:=\operatorname{refine}(\T_{\ell},\T_{\ell})$. Then, $\T_{\ell+1}$ is coarser than $\Tunif$, i.e., 
	\begin{align*}
	\forall T'\in\Tunif\ \exists T\in\T_{\ell+1}:T'\subseteq T.
	\end{align*}
\end{lemma}
\begin{fatproof}
	According to Remark \ref{setting:samesons}, marked elements $T\in\T_{\ell}$ have the same sons in $\T_{\ell+1}$ and $\Tunif$. By definition, and since $\Tunif$ is a conforming triangulation, $\T_{\ell+1}$ is coarser than $\Tunif$.
\end{fatproof}

For a triangulation $\T_{\ell}$, the symbol $\operatorname{refine}(\T_{\ell})$ shall denote the (infinite) set of all triangulations that can be obtained by refining $\T_{\ell}$ finitely many times. Henceforth, we assume that an initial conforming triangulation $\T_0$ is given, and we write $\mathbb{T}:=\operatorname{refine}(\T_0)$ for the set of all possible refinements of $\T_0$.

\subsection{Two-level error estimation and adaptive FEM}

For each node $z\in\SN_{\T_{\ell}}^*$, there is a unique function $\varphi_z^{\ell}\in\S_0^ 1(\T_{\ell})$ that satisfies $\varphi_z^{\ell}(z') = \delta_{zz'}$ for all $z'\in\S_0^1(\T_{\ell})$, where $\delta_{zz'}$ denotes the Kronecker-Delta. The family $\{\varphi_z^{\ell}\ |\ z\in\SN^ *_{\T_{\ell}}\}$ forms a basis of $\S_0^1(\T_{\ell})$ and we call this basis the nodal basis associated to $\T_{\ell}$. Furthermore, for the uniform refinement $\Tunif:=\operatorname{refine}(\T_{\ell},\T_{\ell})$, let $\left\{\varphi^{\ell,u}_z\ |\ z\in\SN_{\Tunif}^*\right\}$ be the nodal basis associated to $\Tunif$. 

For any triangle $T\in\T_{\ell}$, the local contribution of the two-level error estimator (see, e.g., \cite{PRS20}) is defined as
\begin{align}\label{setting:estimator}
	\tau_{\ell}(T)^2:=\sum_{z\in\SN_{\Tunif}^*(T)\setminus\SN_{\T_{\ell}}}\tau_{\ell}\left(\varphi^{\ell,u}_z\right)^2, \text{ where } \tau_{\ell}\left(\varphi^{\ell,u}_z\right) := \frac{\left|\SCP{f}{\varphi^{\ell,u}_z}{L^2(\Omega)}-a\left(u_{\ell},\varphi^{\ell,u}_z\right)\right|}{\vertiii{\varphi^{\ell,u}_z}},
\end{align}

and where $\SN_{\Tunif}^*(T)$ is the set of all vertices of sons of $T$ that do not lie on the boundary of $\Omega$.
Furthermore, for any subset $\mathcal{U}_{\ell}\subseteq \T_{\ell}$, we set $\tau_{\ell}^2(\mathcal{U}_{\ell}):= \sum_{T\in\mathcal{U}_{\ell}}\tau_{\ell}(T)^2$ and $\tau_{\ell}:=\tau_{\ell}(\T_{\ell})$. 

\begin{remark}\label{setting:Galerkinremark}
	In \cite{AFFKP15}, it was noted that there holds $\tau_{\ell}(\varphi^{\ell,u}_z) = \vertiii{\mathbb{G}_z(u_{\ell,u}-u_{\ell})}$, where $\mathbb{G}_z:\S_0^ 1(\Tunif)\rightarrow \operatorname{span}(\varphi^{\ell,u}_z)$ is the Galerkin projection onto $\operatorname{span}(\varphi^{\ell,u}_z)$ and $u_{\ell,u}\in\S_0^ 1(\Tunif)$ is the Galerkin solution corresponding to the uniform refinement $\Tunif$.
\end{remark}

We employ an adaptive algorithm of the form {\texttt{SOLVE}--\texttt{ESTIMATE}--\texttt{MARK}--\texttt{REFINE}}. As error estimators we use the two-level indicators defined above. The marking step is done by using Dörfler marking \cite{Dorfler}, which in fact can be done in linear complexity, \cite{PP20}. Finally, mesh refinement is done with newest vertex bisection, \cite{Stevenson2,NVB}.	

	\begin{algo}\label{SEMR}
		
		 Input:  Conforming triangulation $\T_{0}$, right-hand side $f$, refinement parameter $\theta\in (0,1]$.

		For {$\ell = 0,1,2,...$} do
		{\begin{enumerate}
				\item Compute the unique solution $u_{\ell}\in\S_0^1(\T_{\ell})$ of \eqref{setting:discrweakform}.
				\item For each $T\in\T_{\ell}$, compute the two-level estimator contribution $\tau_{\ell}(T)$ defined in \eqref{setting:estimator} and the total error estimator $\tau_{\ell}$.
				\item 
				Find a set $\mathcal{M}_{\ell}\subseteq \T_{\ell}$ of minimal cardinality such that
					$\tau_{\ell}(\mathcal{M}_{\ell})^2\geq \theta \tau_{\ell}^2$.
				\item Generate the refined mesh $\T_{\ell+1} := \operatorname{refine}(\T_{\ell},\mathcal{M}_{\ell})$.		
				
		\end{enumerate}}
	Output: Sequence of triangulations $(\T_{\ell})_{\ell\in\N_0}$ with associated discrete solutions $(u_{\ell})_{\ell\in\N_0}$.
		
	\end{algo}
	
Regarding convergence of adaptive algorithms with optimal algebraic rates for model problems with integer order as well as for boundary element discretizations of integral equations, we refer to \cite{BDD04,Stevenson,casconfund,CFPP14,Gantumur17} and to \cite{FMP21a} for the fractional Laplacian.

\subsection{Saturation assumptions}

As in \cite{PRS20}, the validity of our main results hinges on certain saturation assumptions. We say that a sequence $(\T_{\ell})_{\ell\in\N_0}$ of triangulations satisfies the \textbf{weak saturation assumption}, if there exists a $q\in (0,1)$ such that  
\begin{align}\label{setting:weaksaturation}
	\vertiii{u-u_{\ell,u}}\leq q\vertiii{u-u_{\ell}}
\end{align} for all $\ell\in\N_0$. Here, $u_{\ell,u}\in\S_0^1(\Tunif)$ denotes the Galerkin solution corresponding to the uniform refinement $\Tunif:=\operatorname{refine}(\T_{\ell}, \T_{\ell})$. 

Furthermore, we say that a sequence $(\T_{\ell})_{\ell\in\N_0}$ of triangulations satisfies the \textbf{strong saturation assumption}, if there exist constants $0<\kappa\leq q <1$ such that for all $\ell\in\N_0$ and all $\T_h\in\operatorname{refine}(\T_{\ell})$ satisfying $\vertiii{u-u_h}\leq \kappa \vertiii{u-u_{\ell}}$, there exists a set $\mathcal{M}_H\subseteq \T_{\ell}\setminus \T_h$ such that $\T_H:=\operatorname{refine}(\T_{\ell},\mathcal{M}_H)$ satisfies 
\begin{align}\label{setting:strongsaturation}
	\vertiii{u-u_H}\leq q\vertiii{u-u_{\ell}} \text{ and } \T_{\ell}\setminus \T_H\subseteq \T_{\ell}\setminus \T_h.
\end{align}
\begin{remark}
	In \cite{PRS20}, the strong saturation assumption was defined differently. There, the strong saturation assumption reads as: Suppose that there exist constants $0<\kappa\leq q<1$ such that for all $\ell\in\N_0$ and $\T_h\in\operatorname{refine}(\T_{\ell})$, the mesh $\T_H:=\operatorname{refine}(\T_{\ell}, \T_{\ell}\setminus\T_h)$ satisfies that 
	\begin{align}\label{strongsatpraet}
		\vertiii{u-u_h}\leq \kappa \vertiii{u-u_{\ell}} \text{ implies } \vertiii{u-u_H}\leq q \vertiii{u-u_{\ell}}.
	\end{align}
	If $\T_{\ell+1}:=\operatorname{refine}(\T_{\ell}, \mathcal{M}_{\ell})$ denotes the coarsest NVB refinement, such that every marked element $T\in\mathcal{M}_{\ell}$ is bisected at least once, both strong saturation assumptions are in fact equivalent, since then, the mesh $\T_H:=\operatorname{refine}(\T_{\ell}, \T_{\ell}\setminus\T_h)$ satisfies $\T_{\ell}\setminus \T_H\subseteq \T_{\ell}\setminus \T_h$. If $\T_{\ell+1}:=\operatorname{refine}(\T_{\ell}, \mathcal{M}_{\ell})$ is defined as in Section \ref{pattern}, however, assumption \eqref{setting:strongsaturation} is stronger than \eqref{strongsatpraet}.
\end{remark}

\subsection{Main results}

Our first main result states that discrete efficiency, discrete reliability and stability holds for the two-level error indicator.

\begin{theorem}\label{thm:mainresult1}
	 Let $\T_{\ell}\in\mathbb{T}$ and $\T_{\ell+1}:=\operatorname{refine}(\T_{\ell},\mathcal{M}_{\ell})$ for some subset $\mathcal{M}_{\ell}\subseteq \T_{\ell}$, as well as $\T_m\in\operatorname{refine}(\T_{\ell})$. Let $u_{\ell}$, $u_{\ell+1}$ and $u_m$ denote the Galerkin solutions of \eqref{setting:discrweakform} corresponding to $\T_{\ell}$, $\T_{\ell+1}$ and $\T_m$, respectively. Then,
	\begin{itemize}
		\item[(a)] there exists a constant $C_{\text{eff}}>0$ such that
	\begin{align*}
		\tau_{\ell}(\mathcal{M}_{\ell})\leq C_{\text{eff}}\  \vertiii{u_{\ell+1}-u_{\ell}},
	\end{align*}
   \item[(b)] there exists a constant $C_{\text{rel}}>0$ such that
   \begin{align*}
   		\vertiii{u_{\ell+1}-u_{\ell}}\leq C_{\text{rel}}\  \tau_{\ell}(\T_{\ell}\setminus\T_{\ell+1}).
   \end{align*}
   \item[(c)] Assume $d=2$. Then, there exists a constant $C_{\text{stab}}>0$ such that
   \begin{align*}
   		\left|\tau_{\ell}(\T_m\cap \T_{\ell})-\tau_m(\T_m\cap\T_{\ell})\right|\leq C_{\text{stab}}\ \vertiii{u_m-u_{\ell}}.
   \end{align*}
\end{itemize}
	The constants $C_{\text{eff}}$, $C_{\text{rel}}$ and $C_{\text{stab}}$ depend only on $\Omega$, $s$, $d$ and the $\gamma$-shape regularity of $\T_{\ell}$.
\end{theorem}

\begin{remark}
	The triangle inequality shows $\vertiii{u_{\ell+1}-u_{\ell}}\leq 2\vertiii{u_{\ell}-u}$ and employment of Dörfler marking means $\theta\tau_{\ell}^2\leq \tau_{\ell}\left({\mathcal{M}_{\ell}}\right)^2$. That is, Theorem~\ref{thm:mainresult1} (a) together with Dörfler marking leads to efficiency of the two-level error indicator.
	Furthermore, the triangle inequality and the weak saturation assumption prove $(1-q)\vertiii{u-u_{\ell}}\leq\vertiii{u_{\ell,u}-u_{\ell}}$. In combination with Theorem~\ref{thm:mainresult1} (b), this shows that the two-level error estimator is reliable, provided that the weak saturation assumption holds.
\end{remark}

\begin{remark}
	In order to prove Theorem \ref{thm:mainresult1} (a) it will be essential to use the refinement strategy $\operatorname{refine}_{\operatorname{FSW}}(\cdot)$, since this refinement strategy guarantees that the nodal basis functions corresponding to newly introduced nodes on marked elements are the same in $\S_0^1(\T_{\ell+1})$ and $\S_0^1(\Tunif)$. The proof of Theorem \ref{thm:mainresult1} (b) stays valid, even if one uses the refinement strategy $\operatorname{refine}_{\operatorname{EGP}}(\cdot)$ instead.
\end{remark}

In combination with the general frameworks presented in \cite{PRS20, CFPP14} and the weak saturation assumption, Theorem \ref{thm:mainresult1} guarantees linear convergence of Algorithm \ref{SEMR}.

\begin{theorem}\label{thm:mainresult2}
	Assume that the output $(\T_{\ell})_{\ell\in\N_0}$ of Algorithm \ref{SEMR} satisfies the weak saturation assumption \eqref{setting:weaksaturation}.
	Then, the sequence of Galerkin solutions $(u_{\ell})_{\ell\in\N_0}$ associated to $(\T_{\ell})_{\ell\in\N_0}$ converges linearly, i.e., there exists a constant $\kappa \in (0,1)$, such that 
	\begin{align*}
		\vertiii{u-u_{\ell+1}}\leq \kappa \vertiii{u-u_{\ell}}
	\end{align*}
	for all $\ell\in\N_0$. The constant $\kappa$ depends only on $\theta$, $C_{\text{eff}}$, $C_{\text{rel}}$ and the constant $q$ from the weak saturation assumption.
\end{theorem}
\begin{fatproof}
	Follows from Theorem \ref{thm:mainresult1} (a), (b) and \cite[Theorem 2.6]{PRS20}.
\end{fatproof}

Furthermore, in two dimensions, the frameworks from \cite{PRS20, CFPP14}, Theorem \ref{thm:mainresult1} and the strong saturation assumption lead to convergence with optimal algebraic convergence rates.


\begin{theorem}\label{thm:mainresult3}
 Assume $d=2$. Consider the output $(\T_{\ell})_{\ell\in\N_0}$ of Algorithm \ref{SEMR} for a fixed refinement parameter $\theta\in (0,1]$, and assume that the strong saturation assumption \eqref{setting:strongsaturation} holds.
	
	Then, there exists $\theta_{\text{opt}}\in (0,1]$, such that, if $\theta<\theta_{\text{opt}}$, the sequence $(u_{\ell})_{\ell\in\N_0}$ of discrete solutions associated with $(\T_{\ell})_{\ell\in\N_0}$ converges with the best possible algebraic rate. That is, for every $t>0$, there exists a constant $C_{\text{opt}}>0$ depending only on $t$, $\theta$, $\T_0$, $C_{\text{eff}}$, $C_{\text{rel}}$, $C_{\text{stab}}$, and the constants $\kappa$ and $q$ from the strong saturation assumption, such that 
	\begin{equation*}
		C_{\text{opt}}^{-1}\norm{u}{\mathbb{A}_t} \leq \sup_{\ell\in\N_0}\left(\#\T_{\ell}\right)^t\vertiii{u-u_{\ell}} \leq C_{\text{opt}}\norm{u}{\mathbb{A}_t},
	\end{equation*}
	where $\norm{u}{\mathbb{A}_t}$ is the so-called approximation constant (\cite{casconfund}) of $u$, 
	\begin{align*}
		\norm{u}{\mathbb{A}_t}:= \sup_{N\in\N_0} (N+1)^t \min_{\{\T_{\text{opt}}\in\mathbb{T}\ |\  \#\T_{\text{opt}}-\#\T_0 \leq N\}} \vertiii{u-u_{\T_{\text{opt}}}}.
	\end{align*}
\end{theorem}
\begin{fatproof}
	Follows from Theorem \ref{thm:mainresult1} (a)-(c) and \cite[Theorem 2.9]{PRS20}.
\end{fatproof}

	\section{Proof of Theorem \ref{thm:mainresult1}}

The proof of our main theorem requires some auxiliary results. As was pointed out in \cite{FMP21a}, one major difficulty in the discussion of a posteriori error estimators for the fractional Laplacian is, that, for $w_{\ell}\in\S_0^1(\T_{\ell})$ and $3/4\leq s<1$, the expression $(-\Delta)^sw_{\ell}$ is in general no longer in $L^2(\Omega)$. In order to overcome this difficulty, a certain weight-function is introduced. For a triangulation $\T_{\ell}$ of $\Omega$ and $s\in (0,1)$ the weight-function $\widetilde{h}_{\ell}^s$ is defined as
\begin{align*}
\widetilde{h}^s_{\ell} := \left\{\begin{array}{ll} h_{\ell}^s & \text{ for } 0<s\leq 1/2 \\
h_{\ell}^{1/2}\omega_{\ell}^{s-1/2} & \text{ for }1/2<s<1,\end{array}\right. 
\end{align*}
where 
\begin{align*}
	\omega_{\ell}(x) := \inf_{T\in\T_{\ell}}\inf_{y\in\partial T} |x-y|, \hskip 0.3cm x\in\Omega
\end{align*}
is the distance from the mesh skeleton. 

\subsection{Properties of the weight-function}

The following properties of the weight-function have been proved in \cite{FMP21a}.

\begin{theorem}\label{inversestimate}
	Consider a triangulation $\T_{\ell}$ of a bounded Lipschitz domain $\Omega\subseteq \R^ d$ and $s\in (0,1)$. 
	\begin{itemize}
		\item[(a)] For all $w_{\ell}\in\S_0^1(\T_{\ell})$, there holds $\widetilde{h}_{\ell}^s(-\Delta)^sw_{\ell}\in L^2(\Omega)$.
		\item[(b)] There exists a constant $C>0$ depending only on $\Omega$, $s$, $d$, and the $\gamma$-shape regularity of $\T_{\ell}$ such that 
		\begin{align*}
			\left\|\widetilde{h}_{\ell}^s(-\Delta)^sw_{\ell}\right\|_{L^2(\Omega)}\leq C\norm{w_{\ell}}{\exspace{\Omega}}
		\end{align*}
		for all $w_{\ell}\in\S_0^1(\T_{\ell})$, i.e., an inverse estimate is available for the fractional Laplacian.
	\end{itemize}
\end{theorem}

 For $T\in\T_{\ell}$ with vertices $z_0, z_1, z_2$ (or $z_0, z_1, z_2, z_3$ in the case $d=3$), let $F_T:\TR\rightarrow T$ be the unique diffeomorphism that satisfies $F(e_i) = F(z_i)$, where $e_i$ shall be the zero vector for $i=0$ and the $i$-th canonical basis vector otherwise. Without loss of generality we may assume that the ordering of the vertices of $T$ is such that the edge $\overline{z_0z_1}$ corresponds to the refinement edge of $T$. That is, $F_T$ maps $\overline{e_0e_1}$ to the refinement edge of $T$.
 
 \medskip 
 
 For the reference element $\TR$, we define
 \begin{align*}
 \omegaref(x) := \inf_{y\in\partial \TR}|x-y|, \hskip 0.3cm x\in\TR.
 \end{align*}
  
 According to \cite[Proof of Lemma 3.2]{FMP21a}, there holds the pointwise estimate
 \begin{align}\label{proofs:scalingfundamental}
 	h_{\ell}(T)\omegaref\simeq \omega_{\ell}\circ F_T,
 \end{align}
 where the hidden constants depend only on $d$ and the $\gamma$-shape regularity of $\T_{\ell}$.

 Consider a triangulation $\T_{\ell}$ with uniform refinement $\Tunif:=\operatorname{refine}(\T_{\ell},\T_{\ell})$. For an element $T\in\T_{\ell}$, the refinement pattern $\Pat$ of $T$ is the set 
 \begin{align*}
 	\Pat := \left\{F_T^ {-1}(T_1),...,F_T^ {-1}(T_N)\right\},
 \end{align*}
 where $T_1,...,T_N$ are the sons of $T$ in $\Tunif$ and $N=4$ in 2D and $N=18$ or $N=20$ (depending on the refinement type of $T$) in 3D.

 For the refinement pattern of an element $T$, we define 
 \begin{align*}
 \omega_{\Pat}(x) := \inf_{T'\in\Pat}\inf_{y\in\partial T'}|x-y|, \hskip 0.3cm x\in\TR.
 \end{align*}
If $d=2$, every element $T$ has the same refinement pattern, i.e., the refinement pattern is independent of $T$. In the three-dimensional case, each element has an associated refinement type $\beta\in\{0,1,2\}$, and the refinement pattern of an element $T$ only depends on its refinement type, i.e., $T$ has one of three possible refinement patterns.
 
   Consider a triangulation $\T_{\ell}$ with uniform refinement $\Tunif:=\operatorname{refine}(\T_{\ell},\T_{\ell})$. Using \eqref{setting:geometric}, arguing as in \cite[Proof of Lemma 3.2]{FMP21a} leads to 
 \begin{align}\label{proofs:scalingfundamental2}
 	(h_{\ell,u}\circ F_T)\omega_{\Pat}\simeq \omega_{\ell,u}\circ F_T
 \end{align}
 for all $T\in\T_{\ell}$, where $h_{\ell,u}$ and $\omega_{\ell,u}$ denote the mesh-width and distance-to-skeleton function associated to $\Tunif$. The hidden constants in \eqref{proofs:scalingfundamental2} depend only on $d$ and the $\gamma$-shape regularity of $\T_{\ell}$.

 \begin{lemma}\label{proofs:hequivalence}
 	Let $s\in(0,1)$ and $\T_{\ell}$ be a triangulation of $\Omega$ and $\T_{\ell+1} := \operatorname{refine}(\T_{\ell}, \mathcal{M}_{\ell})$ for some subset $\mathcal{M}_{\ell}\subseteq \T_{\ell}$. Then, for all $v_{\ell+1}\in \S_0^1(\T_{\ell+1})$ and all $T\in\T_{\ell}$, there holds 
 	\begin{align}\label{proofs:weightequivalence}
 		\left\|\widetilde{h}_{\ell+1}^{-s}v_{\ell+1}\right\|_{L^ 2(T)}\simeq \left\|\widetilde{h}_{\ell}^{-s}v_{\ell+1}\right\|_{L^ 2(T)},
 	\end{align}
 	where the hidden constants depend only on $s$, $d$ and the $\gamma$-shape regularity of $\T_{\ell}$.
 \end{lemma}

\begin{fatproof}
	Due to \eqref{setting:geometric}, the statement is clear for $s\leq1/2$. 
	Therefore, consider $s>1/2$ in the following. Again, due to \eqref{setting:geometric}, the inequalities \eqref{proofs:weightequivalence} are equivalent to 
	\begin{align}\label{proofs:tmpestimate}
	\norm{\omega_{\ell+1}^{1/2-s}v_{\ell+1}}{L^2(T)}\simeq \norm{\omega_{\ell}^{1/2-s}v_{\ell+1}}{L^2(T)}.
	\end{align}
	Due to the pointwise inequality $\omega_{\ell+1}\leq\omega_{\ell}$ and $s>1/2$, bounding the right-hand side by the left-hand side in \eqref{proofs:tmpestimate} is clear and we only have to prove the converse estimate. 
	
{\bf Step~1:}	The first step is to prove the lower bound in \eqref{proofs:tmpestimate} under the assumption that $\mathcal{M}_{\ell} = \T_{\ell}$, i.e., that $\T_{\ell+1} = \Tunif$. 
	
	To this end, let $T\in\T_{\ell}$ be given and consider the induced distance-function $\omega_{\Pat}$ on the reference element. Furthermore, we define the space 
 \begin{align*}
 \S^1(\Pat) := \left\{u\in C(\TR)\ |\ u\restrict{T}\in \mathcal{P}^1(T)\text{ for all }T\in\Pat\right\}.
 \end{align*}

  For $v_{\text{ref}}\in\S^ 1(\Pat)$ consider the weighted $L^ 2$-norms 
	\begin{align*}
		\norm{v_{\text{ref}}}{\omega_{\Pat}}^2 := \int_{\TR}(\omega_{\Pat}(x))^ {1-2s}v_{\text{ref}}^ 2\ \diff x,\ \text{ and } 	\norm{v_{\text{ref}}}{\omegaref}^2 := \int_{\TR}(\omegaref)^ {1-2s}v_{\text{ref}}^ 2\ \diff x.
	\end{align*}
	Since $\S^1(\Pat)$ is finite-dimensional and $\Pat$ is independent of $T$ (when $d=2$) or depends only on its refinement type $\beta\in\{0,1,2\}$ (when $d=3$), the norms $\norm{\cdot}{\omega_{\Pat}}$ and $\norm{\cdot}{\omegaref}$ are pairwise equivalent and there exists a constant $C>
	0$ depending only on $s$ and $d$ such that $C^{-1}\norm{\cdot}{\omega_{\Pat}}\leq\norm{\cdot}{\omegaref}\leq C\norm{\cdot}{\omega_{\Pat}}$.
	In combination with \eqref{proofs:scalingfundamental} and \eqref{proofs:scalingfundamental2}, a scaling argument proves the lower bound in \eqref{proofs:tmpestimate}, provided that $\T_{\ell+1}$ is the uniform refinement of $\T_{\ell}$.
	
{\bf Step~2:}	In the more general case, consider the distance-to-skeleton function $\omega_{\ell,u}$, which is induced by the uniform refinement $\Tunif$. From Lemma \ref{setting:uniforminref} we infer the pointwise inequality $\omega_{\ell,u}\leq\omega_{\ell+1}$, which, together with $s>1/2$ and the first step of the proof therefore implies 
	\begin{align*}
		\norm{\omega_{\ell+1}^{1/2-s}v_{\ell+1}}{L^2(T)}\leq \norm{\omega_{\ell,u}^{1/2-s}v_{\ell+1}}{L^2(T)}\stackrel{step~1}{\lesssim} \norm{\omega_{\ell}^{1/2-s}v_{\ell+1}}{L^2(T)}
	\end{align*}
	for all $v_{\ell+1}\in \S_0^1(\T_{\ell+1})$, which concludes the proof.
	
\end{fatproof}

\subsection{Scott-Zhang projection}

We recall that, for a triangulation $\T_{\ell}$, the Scott-Zhang projection $J_{\ell}: \exspace{\Omega}\rightarrow \S_0^1(\T_{\ell})$ is defined as follows: For any node $z\in\mathcal{N}_{\T_{\ell}}^*$, let $T_z\in\T_{\ell}$ be an element that has $z$ as vertex. Furthermore, let $\varphi^{\ell}_z\in\S_0^1(\T_{\ell})$ denote the nodal basis function associated to $z$ and let $\varphi_z^{\ell*}\in\mathcal{P}^1(T_z)$ be the unique dual function that satisfies $\SCP{\varphi_z^{\ell*}}{\varphi^{\ell}_{z'}}{L^2(T_z)}=\delta_{zz'}$ for all $z'\in\mathcal{N}_{\T_{\ell}}^*(T_z)$, where $\delta_{zz'}$ is the Kronecker-Delta. Then, for $v\in\exspace{\Omega}$, 
\begin{align*}
	J_{\ell}(v) := \sum_{z\in\mathcal{N}_{\T_{\ell}}^*} \left(\int_{T_z}\varphi_z^{\ell*}v\ \diff x\right) \varphi^{\ell}_z.
\end{align*}

The Scott-Zhang projection has the following stability and approximation properties (see \cite[Lemma 3.2]{FMP21a}): 

\begin{theorem}\label{scproperties}
	For $s\in (0,1)$, there exists a constant $C>0$ depending only on $s$, the $\gamma$-shape regularity of $\T_{\ell}$, $\Omega$ and $d$, such that the following statements hold:
	\begin{enumerate}
		\item[(a)] $\norm{J_{\ell}v}{\widetilde{H}^s(\Omega)}\leq C\norm{v}{\widetilde{H}^s(\Omega)}$ for all $v\in \widetilde{H}^s(\Omega)$.
		\item[(b)]$\norm{\widetilde{h}_{\ell}^{-s}(1-J_{\ell})v}{L^2(\Omega)}\leq C\norm{v}{\widetilde{H}^s(\Omega)}$ for all $v\in\widetilde{H}^s(\Omega)$.
	\end{enumerate}
\end{theorem}

For a triangulation $\T_{\ell}$, let $I_{\ell}: C(\Omega)\rightarrow \S_0^1(\T_{\ell})$ be the nodal interpolation operator, i.e., $I_{\ell}(v) := \sum_{z\in\mathcal{N}_{\T_{\ell}}^*}v(z)\varphi^{\ell}_z$. The next theorem shows that, on the discrete space $\S_0^1(\T_{\ell})$, the nodal interpolation operator and the Scott-Zhang projection are in some sense equivalent.

\begin{theorem}\label{thm:NodalSZ}
	Let $\T_{\ell}\in\mathbb{T}$ with uniform refinement $\Tunif:=\operatorname{refine}(\T_{\ell}, \T_{\ell})$. Then, for all $v_{\ell,u}\in\S_0^1(\Tunif)$, there holds the estimate
	\begin{equation*}
		\norm{\widetilde{h}_{\ell}^{-s}(1-I_{\ell})v_{\ell,u}}{L^2(\Omega)}\simeq \norm{\widetilde{h}_{\ell}^{-s}(1-J_{\ell})v_{\ell,u}}{L^2(\Omega)},
	\end{equation*}
	where the hidden constants may depend only on $s$, $d$, and the $\gamma$-shape regularity of $\T_{\ell}$.

\end{theorem}
\begin{fatproof}
	 It is enough to prove the element-wise estimate
	\begin{equation*}
	\norm{\widetilde{h}_{\ell}^{-s}(1-I_{\ell})v_{\ell,u}}{L^2(T)}\simeq \norm{\widetilde{h}_{\ell}^{-s}(1-J_{\ell})v_{\ell,u}}{L^2(T)}
	\end{equation*}
	for all $T\in\T_{\ell}$.
	The proof of this estimate is split into four steps.
	
	\bigskip
	
		\textbf{Step 1}. 
		For any $T\in\T_{\ell}$, we have already considered the affine diffeomorphism $F_T:\TR\rightarrow T$. Let $V_T:\R^d\rightarrow\R^d$ be the unique affine extension of $F_T$ to $\R^d$. Note that $F_T$ and $V_T$ have the same Jacobi matrix. For $T\in\T_{\ell}$, we define 
		\begin{align*}
			\Pi_{\ell}^{\text{ref}}[T] := \{V_T^{-1}(T')\ |\ \overline{T}\cap\overline{T'}\neq\emptyset\},
		\end{align*}
		as well as 
		\begin{align*}
		\Omega_{\ell}^{\text{ref}}[T] := \operatorname*{interior}\bigcup\Big\{\overline{V_T^{-1}(T')}\ |\ \overline{T}\cap\overline{T'}\neq\emptyset\Big\}.
		\end{align*}
		The affine mapping $P_T:=V_T\vert_{\Omega_{\ell}^{\text{ref}}[T]}$ then maps $\Omega_{\ell}^{\text{ref}}[T]$ onto $\Omega_{\ell}[T]$, and we may consider the space
		\begin{align*}
			\S^1(P_T^{-1}(\Tunif)) := \left\{v_{\ell,u}\circ P_T\ \big|\ v_{\ell,u}\in\S_0^1(\Tunif)\right\}.
		\end{align*}
		This definition suggests that $\S^1(P_T^{-1}(\Tunif))$ depends on the element $T$. A closer look, however, reveals that it depends only on the shape of the element patch of $T$ and the distribution of the refinement edges of elements in the patch (and the refinement types $\beta\in\{0,1,2\}$ of the elements in the patch in $d=3$). The amount of possible shapes of element patches can be controlled in terms of $\gamma$-shape regularity. This shows, that $\S^1(P_T^{-1}(\Tunif))$ is one of only $M$ possible spaces, where $M$ depends on $\T_0$, $d$ and $\gamma$-shape regularity, but is independent of the number of elements in $\T_{\ell}$.


\bigskip 

\textbf{Step 2.} Let an element $T\in\T_{\ell}$ with associated function space $\S^1(P_T^{-1}(\Tunif))$ be given. We define an operator
	 $I^{\text{ref}}: \S^1(P_T^{-1}(\Tunif))\rightarrow \S^1(P_T^{-1}(\Tunif))$ by
	\begin{align*}
	I^{\text{ref}}(v_{\ell,u}\circ P_T) := \sum_{z\in\mathcal{N}_{\T_{\ell}}^*(T)}v_{\ell,u}(z)\varphi^{\ell}_{z}\circ P_T.
	\end{align*}
	
	
	We observe that for all $v_{\ell,u}\in\S_0^1(\Tunif)$ there holds $((1-I_{\ell})v_{\ell,u})\circ F_T = (1-I^{\text{ref}})(v_{\ell,u}\circ P_T)\vert_{\TR}$.
	
	Furthermore, for $T\in\T_{\ell}$ with associated function space $\S^1(P_T^{-1}(\Tunif))$, we define an operator $J_{T}^{\text{ref}}: \S^1(P_T^{-1}(\Tunif))\rightarrow \S^1(P_T^{-1}(\Tunif))$ by
	\begin{align*}
		J_T^{\text{ref}}(v_{\ell,u}\circ P_T) := \sum_{z\in\mathcal{N}_{\T_{\ell}}^*(T)} \left(\int_{T_z}\varphi_z^{\ell*}v_{\ell,u}\ \diff x\right) \varphi^{\ell}_z\circ P_T,
	\end{align*}
	
	where the elements $T_z$ are the averaging elements associated to vertices $z$ of $T$ chosen in the definition of $J_{\ell}$.
	We keep the subscript $T$ here, since this operator is also dependent on the choice of the averaging elements $T_z$. In other words: Not every information needed to define $J_T^{\text{ref}}$ is incorporated in the function space $\S^1(P_T^{-1}(\Tunif))$, one still has an option to choose different averaging elements. This is in contrast to the operator $I^{\text{ref}}$, which is determined by $\S^1(P_T^{-1}(\Tunif))$.
	
	By definition, $J_T^{\text{ref}}$ satisfies $((1-J_{\ell})v_{\ell,u})\circ F_T = (1-J_T^{\text{ref}})(v_{\ell,u}\circ P_T)\vert_{\TR}$ for all functions $v_{\ell,u}\in\S_0^1(\Tunif)$. 
	
	\bigskip 
	
	\textbf{Step 3.} For $T\in\T_{\ell}$, the functions
	\begin{equation*}
		\norm{g}{I} := \norm{\omega_{\text{ref}}^{1/2-s}(1-I^{\text{ref}})g}{L^2(\TR)}
	\end{equation*}
	and
	\begin{equation*}
		\norm{g}{J, T}:= \norm{\omega_{\text{ref}}^{1/2-s}(1-J_T^{\text{ref}})g}{L^2(\TR)}
	\end{equation*}
	are seminorms on the associated space $\S^1(P_T^{-1}(\Tunif))$. For each $T\in\T_{\ell}$, the kernels of $\norm{\cdot}{I}$ and $\norm{\cdot}{J, T}$ coincide, and, since the number of elements $T$ in $\T_{\ell}$ is finite, there exists a constant $C_{\T_{\ell},s}>0$ such that
	\begin{equation*}
		C_{\T_{\ell},s}^{-1}\norm{g}{J, T}\leq \norm{g}{I} \leq C_{\T_{\ell},s}\norm{g}{J, T}
	\end{equation*}
	for all $g\in\S^1(P_T^{-1}(\Tunif))$ and all $T\in\T_{\ell}$. The aim of the third step is to show that there is a constant $C_s>0$, which may depend on $s$, $\gamma$-shape regularity and $d$ but not on $\T_{\ell}$, such that
	\begin{equation}\label{Tindependance}
		C_{s}^{-1}\norm{g}{J, T}\leq \norm{g}{I} \leq C_{s}\norm{g}{J, T}
	\end{equation}
	for all $g\in\S^1(P_T^{-1}(\Tunif))$.
	
	To this end, let an element $T\in\T_{\ell}$ with associated space $\S^1(P_T^{-1}(\Tunif))$ be given.
	The set $\left\{\varphi_{z'}^{\ell,u}\circ P_T\ |\ z'\in\mathcal{N}_{\Tunif}^*\right\}$ is a spanning set of $\S^1(P_T^{-1}(\Tunif))$. Note the amount of nonzero functions in this spanning set is bounded by a number depending only on $d$ and $\gamma$-shape regularity.
	For each function in this set, the definition of $J_T^{\text{ref}}$ shows
	\begin{align*}
	J_T^{\text{ref}}(\varphi_{z'}^{\ell,u}\circ P_T) := \sum_{z\in\mathcal{N}_{\T_{\ell}}^*(T)} \left(\int_{T_{z}}\varphi_{z}^{\ell*}\varphi_{z'}^{\ell,u}\ \diff x\right) \varphi^{\ell}_{z}\circ P_T.
	\end{align*}
	
	The function $\varphi_{z}^{\ell}\circ P_T$ is a nodal basis function associated to a node of $\TR$, i.e., it is one of only $d+1$ functions.
	The key is to take a closer look on the coefficients $\left(\int_{T_{z}}\varphi_{z}^{\ell *}\varphi_{z'}^{\ell,u}\ \diff x\right)$. If $\varphi_{z'}^{\ell,u}$ does not have support in $T_z$, then the coefficient is zero. We observe that, if $\varphi_{z'}^{\ell,u}$ has support in $T_z$, then  $\varphi_{z'}^{\ell,u}\circ F_{T_z}$ is one of only $M$ functions, where $M=6$ in the case $d=2$ or $M=14$ for $d=3$. This is due to the fact, that $\varphi_{z'}^{\ell,u}\circ F_{T_z}$ corresponds to a nodal basis function  associated to vertices of $\TR$ or nodes on the midpoints of edges or contained in faces of $\TR$.
	
	The transformation rule shows that $h_{\ell}(T_z)^d(\varphi_z^{\ell *}\circ F_{T_z})\in\mathcal{P}^1(\TR)$ is a dual function on the reference element, and 
	\begin{align*}
		\frac{1}{2}\int_{T_{z}}\varphi_{z}^{\ell *}\varphi_{z'}^{\ell,u}\ \diff x 
		= \int_{\TR}h_{\ell}(T_z)^d\left(\varphi_z^{\ell *}\circ F_{T_z}\right)\ \left(\varphi_{z'}^{\ell,u}\circ F_{T_z}\right)\ \diff x.
	\end{align*}
    Since there are only $d+1$ dual functions in $\mathcal{P}^1(\TR)$, this implies that the number of possible values of the coefficient $\int_{T_{z}}\varphi_{z}^{\ell *}\varphi_{z'}^{\ell,u}\ \diff x $ can be bounded by a constant depending only on $d$.
    
    Since the number of possible coefficients can be bounded independently of $\T_{\ell}$, it follows that the number of functions of the form $J_T^{\text{ref}}(\varphi_{z'}^{\ell,u}\circ P_T)$ is bounded independently of $T$.  Consequently, the number of interpolation operators $J_T^{\text{ref}}$ depends only on $d$ and not on the number of elements in $\T_{\ell}$. In other words: The number of possible operators $J_T^{\text{ref}}$ (and therefore the number of possible seminorms $\norm{\cdot}{J,T}$ on the space $\S^1(P_T^{-1}(\Tunif))$) is bounded independently of the number of elements in $\T_{\ell}$. Together with the fact, that the amount of possible spaces $\S^1(P_T^{-1}(\Tunif))$ is bounded independently of $\T_{\ell}$ (see Step 1),  this shows \eqref{Tindependance}.

	\bigskip
	\textbf{Step 4.} According to step 3, for $T\in\T_{\ell}$ and $g\in \S^1(P_T^{-1}(\Tunif))$, we have  
	\begin{align*}
		\norm{\omega_{\text{ref}}^{1/2-s}(1-I^{\text{ref}})g}{L^2(\TR)} \simeq \norm{\omega_{\text{ref}}^{1/2-s}(1-J_T^{\text{ref}})g}{L^2(\TR)},
	\end{align*}
	where the hidden constants are independent of $T$. In combination with a transformation to the reference element, equation \eqref{proofs:scalingfundamental} and the definition $J_T^{\text{ref}}$ show
	\begin{align*}
		\norm{\widetilde{h}_{\ell}^{-s}(1-I_{\ell})v_{\ell,u}}{L^2(T)}^2 
		&\simeq 2h_{\ell}(T)^{d-2s}\int_{\TR}\omega_{\text{ref}}^{1-2s}((1-I^{\text{ref}})(v_{\ell,u}\circ P_T))^2\ \diff x \\
		&\simeq h_{\ell}(T)^{d-2s}\int_{\TR}\omega_{\text{ref}}^{1-2s}((1-J_T^{\text{ref}})(v_{\ell,u}\circ P_T))^2\ \diff x \\
		&\simeq h_{\ell}(T)^{d-2s}\int_{\TR}\omega_{\text{ref}}^{1-2s}(((1-J_{\ell})v_{\ell,u})\circ F_T)^2\ \diff x.
	\end{align*}
    Transforming back to the reference element (and again using \eqref{proofs:scalingfundamental}) yields 
    \begin{align*}
    	h_{\ell}(T)^{d-2s}\int_{\TR}\omega_{\text{ref}}^{1-2s}(((1-J_{\ell})v_{\ell,u})\circ F_T)^2\ \diff x \simeq \norm{\widetilde{h}_{\ell}^{-s}(1-J_{\ell})v_{\ell,u}}{L^2(T)}^2.
    \end{align*}
	Thus, we have shown
	\begin{align*}
		\norm{\widetilde{h}_{\ell}^{-s}(1-I_{\ell})v_{\ell,u}}{L^2(T)}^2 \simeq \norm{\widetilde{h}_{\ell}^{-s}(1-J_{\ell})v_{\ell,u}}{L^2(T)}^2
	\end{align*}
    for all $T\in\T_{\ell}$ and $v_{\ell,u}\in \S^1(\Tunif)$.
	Summing over $T\in\T_{\ell}$ completes the proof.
\end{fatproof}

Theorem \ref{thm:NodalSZ} yields equivalence of the energy norm and a weighted $L^2$-norm for certain discrete functions.
\begin{lemma}\label{normequivalence}
	Let $\T_{\ell}\in\mathbb{T}$, $\mathcal{M}_{\ell}\subseteq \T_{\ell}$ and $\T_{\ell+1}:=\operatorname{refine}(T_{\ell}, \mathcal{M}_{\ell})$. For any node $z\in \mathcal{N}^*_{\T_{\ell+1}}\setminus\mathcal{N}_{\T_{\ell}}$ consider the associated nodal basis function $\varphi_z^{\ell+1}\in \S_0^1(\T_{\ell+1})$. Then, there holds
	\begin{align*}
	\vertiii{\varphi_z^{\ell+1}}\simeq\norm{\widetilde{h}_{\ell}^{-s}\varphi_z^{\ell+1}}{L^2(\Omega)},
	\end{align*}
	where the hidden constants depend only on $s$, $d$, $\Omega$, and the $\gamma$-shape regularity of $\T_{\ell}$.
\end{lemma}
\begin{fatproof}
	The Cauchy-Schwarz inequality, Lemma \ref{proofs:hequivalence} and Theorem \ref{inversestimate}(b) lead to
	\begin{align*}
	\vertiii{\varphi_z^{\ell+1}}^2 = a\left(\varphi_z^{\ell+1},\varphi_z^{\ell+1}\right)
	&= \int_{\Omega}(-\Delta)^s(\varphi_z^{\ell+1})\widetilde{h}_{\ell+1}^{s} \widetilde{h}_{\ell+1}^{-s} \varphi_z^{\ell+1}\ \diff x \\
	&\leq \norm{\widetilde{h}^s_{\ell+1}(-\Delta)^s(\varphi_z^{\ell+1})}{L^2(\Omega)}\norm{\widetilde{h}_{\ell+1}^{-s}\varphi_z^{\ell+1}}{L^2(\Omega)} \\
	&\lesssim \norm{\widetilde{h}^s_{\ell+1}(-\Delta)^s(\varphi_z^{\ell+1})}{L^2(\Omega)}\norm{\widetilde{h}_{\ell}^{-s}\varphi_z^{\ell+1}}{L^2(\Omega)} \\
	&\lesssim \vertiii{\varphi_z^{\ell+1}}\norm{\widetilde{h}_{\ell}^{-s}\varphi_z^{\ell+1}}{L^2(\Omega)},
	\end{align*}
	which proves the lower inequality. For the reverse inequality note that, for a node $z\in \mathcal{N}^*_{\T_{\ell+1}}\setminus\mathcal{N}_{\T_{\ell}}$, there holds $\varphi_z^{\ell+1}=(1-I_{\ell})\varphi_z^{\ell+1}$. Therefore, Theorem \ref{thm:NodalSZ} and Theorem \ref{scproperties} imply
	\begin{align*}
	\norm{\widetilde{h}_{\ell}^{-s}\varphi_z^{\ell+1}}{L^2(\Omega)} = \norm{\widetilde{h}_{\ell}^{-s}(1-I_{\ell})\varphi_z^{\ell+1}}{L^2(\Omega)} \simeq \norm{\widetilde{h}_{\ell}^{-s}(1-J_{\ell})\varphi_z^{\ell+1}}{L^2(\Omega)}\lesssim \vertiii{\varphi_z^{\ell+1}},
	\end{align*}
	where the last inequality is due to norm equivalence of $\vertiii{\cdot}$ and $\norm{\cdot}{\widetilde{H}^s(\Omega)}$.
\end{fatproof}

\subsection{Proof of Theorem \ref{thm:mainresult1} (a)}

We need a generalization of the definition of the two-level indicators $\tau_{\ell}\left(\varphi_z^{\ell,u}\right)$ defined in \eqref{setting:estimator}. For any refinement $\T_{\ell+1}=\operatorname{refine}(\T_{\ell},\mathcal{M}_{\ell})$ and a nodal basis function $\varphi_z^{\ell+1}\in\S_0^1(\T_{\ell+1})$ associated to a node $z\in\SN_{\T_{\ell+1}}^*(T)\setminus\SN_{\T_{\ell}}$, we set 
\begin{align*}
\tau_{\ell}\left(\varphi^{\ell+1}_z\right) := \frac{\left|\SCP{f}{\varphi^{\ell+1}_z}{L^2(\Omega)}-a\left(u_{\ell},\varphi^{\ell+1}_z\right)\right|}{\vertiii{\varphi^{\ell+1}_z}}.
\end{align*}
The following lemma proves the first part of Theorem \ref{thm:mainresult1}.
\begin{lemma}\label{proofs:localestimate}
	 Let $\T_{\ell}\in\mathbb{T}$, $\mathcal{M}_{\ell}\subseteq \T_{\ell}$ and $\T_{\ell+1}:=\operatorname{refine}(T_{\ell}, \mathcal{M}_{\ell})$. Assume $u_{\ell}$ and $u_{\ell+1}$ are the Galerkin solutions with respect to the triangulations $\T_{\ell}$ and $\T_{\ell+1}$. Then, there holds
	\begin{align}\label{proofs:tmp1}
	\sum_{z\in \mathcal{N}^*_{\T_{\ell+1}}\setminus\mathcal{N}_{\T_{\ell}} }\tau_{\ell}\left(\varphi_z^{\ell+1}\right)^2\lesssim \vertiii{u_{\ell+1}-u_{\ell}}^2
	\end{align}
	as well as
	\begin{align}\label{proofs:discreff}
		\tau_{\ell}(\mathcal{M}_{\ell})\lesssim \vertiii{u_{\ell+1}-u_{\ell}},
	\end{align}
	where the hidden constants depend only on $s$, $d$, $\Omega$, and the $\gamma$-shape regularity of $\T_{\ell}$.
\end{lemma}

\begin{fatproof}
	Consider a node $z\in \mathcal{N}_{\T_{\ell+1}}^*\setminus\mathcal{N}_{\T_{\ell}}$ and an element $T\in\T_{\ell}$ such that $z\in T$. Note that
	\begin{equation*}
	\begin{split}
	\left|\SCP{f}{\varphi_z^{\ell+1}}{L^2(\Omega)}-a\left(u_{\ell},\varphi_z^{\ell+1}\right)\right| &= \left|a\left(u_{\ell+1}-u_{\ell},\varphi_z^{\ell+1}\right)\right|
	 \\
	&= \left|\int_{\Omega}\widetilde{h}_{\ell+1}^s\left((-\Delta)^s (u_{\ell+1}-u_{\ell})\right)\  \widetilde{h}_{\ell+1}^{-s}\varphi_z^{\ell+1}\ \diff x \right|.
	\end{split}
	\end{equation*}
	Due to $\operatorname{supp}(\varphi_z^{\ell+1})$ being contained in the patch $\Omega_{\ell}[T]$, the Cauchy-Schwarz inequality and the identity $(1-I_{\ell})\varphi_z^{\ell+1} = \varphi_z^{\ell+1}$ show
	\begin{equation*}
	\begin{split}
	&\left|\int_{\Omega}\widetilde{h}_{\ell+1}^s\left((-\Delta)^s (u_{\ell+1}-u_{\ell})\right)\  \widetilde{h}_{\ell+1}^{-s}\varphi_z^{\ell+1}\ \diff x \right|
	 =\left|\int_{\Omega_{\ell}[T]}\widetilde{h}_{\ell+1}^s\left((-\Delta)^s (u_{\ell+1}-u_{\ell})\right)\  \widetilde{h}_{\ell+1}^{-s}\varphi_z^{\ell+1}\ \diff x \right| \\
	&\hskip 4cm \leq \norm{\widetilde{h}_{\ell+1}^s(-\Delta)^s(u_{\ell+1}-u_{\ell})}{L^2(\Omega_{\ell}[T])}\norm{\widetilde{h}_{\ell+1}^{-s}(1-I_{\ell})\varphi_z^{\ell+1}}{L^2(\Omega)}. 
	\end{split}
	\end{equation*}
	Furthermore, Lemma \ref{proofs:hequivalence}, Theorem \ref{thm:NodalSZ} and Theorem \ref{scproperties}(b) imply 
	\begin{align*}
		\norm{\widetilde{h}_{\ell+1}^{-s}(1-I_{\ell})\varphi_z^{\ell+1}}{L^2(\Omega)} \lesssim \norm{\widetilde{h}_{\ell}^{-s}(1-I_{\ell})\varphi_z^{\ell+1}}{L^2(\Omega)}\lesssim \norm{\widetilde{h}_{\ell}^{-s}(1-J_{\ell})\varphi_z^{\ell+1}}{L^2(\Omega)} \lesssim \vertiii{\varphi_z^{\ell+1}}.
	\end{align*}
	Thus, for all $z\in \mathcal{N}_{\T_{\ell+1}}^*\setminus\mathcal{N}_{\T_{\ell}}$, we have
	\begin{equation*}
	\begin{split}
	\tau_{\ell}\left(\varphi_z^{\ell+1}\right) = \frac{\left|\SCP{f}{\varphi_z^{\ell+1}}{L^2(\Omega)}-a\left(u_{\ell},\varphi_z^{\ell+1}\right)\right|}{\vertiii{\varphi_z^{\ell+1}}} 
	\lesssim  \norm{\widetilde{h}_{\ell+1}^s(-\Delta)^s(u_{\ell+1}-u_{\ell})}{L^2(\Omega_{\ell}[T])}.
	\end{split}
	\end{equation*}
	Consequently, summing over all $z\in \mathcal{N}_{\T_{\ell+1}}^*\setminus\mathcal{N}_{\T_{\ell}}$, using $\gamma$-shape regularity and Theorem \ref{inversestimate}(b), we conclude
	\begin{equation*}
	\begin{split}
	\sum_{z\in \mathcal{N}_{\T_{\ell+1}}^*\setminus\mathcal{N}_{\T_{\ell}} }\tau_{\ell}\left(\varphi_z^{\ell+1}\right)^2
	&\lesssim \sum_{z\in \mathcal{N}_{\T_{\ell+1}}^*\setminus\mathcal{N}_{\T_{\ell}} }\norm{\widetilde{h}_{\ell+1}^s(-\Delta)^s(u_{\ell+1}-u_{\ell})}{L^2(\Omega_{\ell}[T])}^2 \\
	&\lesssim \norm{\widetilde{h}_{\ell+1}^s(-\Delta)^s(u_{\ell+1}-u_{\ell})}{L^2(\Omega)}^2 
	\lesssim \vertiii{u_{\ell+1}-u_{\ell}}^2.
	\end{split}
	\end{equation*}
	This completes the proof of \eqref{proofs:tmp1}.
		In order to obtain \eqref{proofs:discreff}, we have to compare the nodal basis function $\varphi_z^{\ell+1}$ and $\varphi_z^{\ell,u}$. According to the refinement strategy used in this work, for any element $T\in\mathcal{M}_{\ell}$, there holds $\mathcal{N}_{\T_{\ell+1}}(T) = \mathcal{N}_{\Tunif}(T)$. 
	
	In the two-dimensional case it was observed in \cite[Proof of Theorem 3.2]{PRS20} that, for $T\in\mathcal{M}_{\ell}$ and $z\in \mathcal{N}^*_{\Tunif}(T)\setminus\mathcal{N}_{\T_{\ell}}$, there holds the equality $\varphi_z^{\ell+1} = \varphi_z^{\ell,u}$. 
	
	This equality does also hold in 3D: Consider $T\in\mathcal{M}_{\ell}$ and $z\in \mathcal{N}^*_{\Tunif}(T)\setminus\mathcal{N}_{\T_{\ell}}$. That is, $z$ is either the midpoint of an edge of $T$ or is contained in a face of $T$. Consequently, the nodal basis function $\varphi_z^{\ell+1}$ only has support in elements $T'\in\T_{\ell}$ which share at least one edge with the marked element $T$, i.e,
	\begin{align*}
		\operatorname{supp}(\varphi_z^{\ell+1})\subseteq\bigcup\left\{\overline{T}\ |\ T\in \Lambda(\mathcal{M}_{\ell})\right\}.
	\end{align*} 
	 Since the sons of elements $T'\in\Lambda(\mathcal{M}_{\ell})$ are the same in $\T_{\ell+1}$ and $\Tunif$ (see Remark \ref{setting:samesons}), we conclude that $\varphi_z^{\ell+1}=\varphi_z^{\ell,u}$.

	 Finally, the equalities $\mathcal{N}_{\T_{\ell+1}}(T) = \mathcal{N}_{\Tunif}(T)$ and $\varphi_z^{\ell+1}=\varphi_z^{\ell,u}$ (for elements $T\in\mathcal{M}_{\ell}$ and nodes $z\in \mathcal{N}^*_{\Tunif}(T)\setminus\mathcal{N}_{\T_{\ell}}$) imply 
	\begin{align*}
	\tau_{\ell}(\mathcal{M}_{\ell})^2 = \sum_{T\in\mathcal{M}_{\ell}}\sum_{z\in \mathcal{N}^*_{\T_{\ell+1}}(T)\setminus\mathcal{N}_{\T_{\ell}}}\tau_{\ell}\left(\varphi_z^{\ell+1}\right)^2 \lesssim \sum_{z\in \mathcal{N}^*_{\T_{\ell+1}}\setminus\mathcal{N}_{\T_{\ell}} }\tau_{\ell}\left(\varphi_z^{\ell+1}\right)^2,
	\end{align*}
	and \eqref{proofs:discreff} follows from \eqref{proofs:tmp1}.
	
\end{fatproof} 

\subsection{Proof of Theorem \ref{thm:mainresult1} (b)}

 The proof of Theorem~\ref{thm:mainresult1} (b) requires an auxiliary lemma.

\begin{lemma}\label{inversetriangle}
	Let $\T_{\ell}\in\mathbb{T}$, $\mathcal{M}_{\ell}\subseteq \T_{\ell}$ and $\T_{\ell+1}:=\operatorname{refine}(T_{\ell}, \mathcal{M}_{\ell})$. Consider a function $v_{\ell+1}\in\S_0^1(\T_{\ell+1})$, i.e.,
	\begin{equation*}
		v_{\ell+1} = \sum_{z\in\mathcal{N}_{\T_{\ell+1}}^*} \lambda_z\varphi_z^{\ell+1}
	\end{equation*}
    for certain coefficients $\lambda_z$.
	Then, there holds
	\begin{equation*}
		\sum_{z\in\mathcal{N}_{\T_{\ell+1}}^*}\norm{\widetilde{h}_{\ell}^{-s}\lambda_z\varphi_z^{\ell+1}}{L^2(\Omega)}^2\simeq \norm{\widetilde{h}_{\ell}^{-s}v_{\ell+1}}{L^2(\Omega)}^2.
	\end{equation*}
	The hidden constants depend only on $s$, $d$, and the $\gamma$-shape regularity of $\T_{\ell}$.
	
\end{lemma}

\begin{fatproof}
	The proof is done similarly as the proof of \cite[Lemma 15]{AFFKP15}.  For any element $T\in\T_{\ell}$ the transformation rule and \eqref{proofs:scalingfundamental} show
	\begin{equation*}
		\begin{split}
			\sum_{z\in\mathcal{N}_{\T_{\ell+1}}^*}\norm{\widetilde{h}_{\ell}^{-s}\lambda_z\varphi_z^{\ell+1}}{L^2(T)}^2 &=
			\sum_{z\in\mathcal{N}_{\T_{\ell+1}}^*(T)}\norm{\widetilde{h}_{\ell}^{-s}\lambda_z\varphi_z^{\ell+1}}{L^2(T)}^2 \\ &=2h_{\ell}(T)^d\sum_{z\in\mathcal{N}_{\T_{\ell+1}}^*(T)}\norm{\left(\widetilde{h}_{\ell}^{-s}\lambda_z\varphi_z^{\ell+1}\right)\circ F_T}{L^2(\TR)}^2 \\
			&\simeq  2h_{\ell}(T)^{d-2s}\sum_{z\in\mathcal{N}_{\T_{\ell+1}}^*(T)}\norm{\omega_{\text{ref}}^{1/2-s}\lambda_z\left(\varphi_z^{\ell+1}\circ F_T\right)}{L^2(\TR)}^2. \\
		\end{split}
	\end{equation*}
     The functions $\varphi_z^{\ell+1}\circ F_T$ are linearly independent, therefore norm equivalence on finite-dimensional spaces shows 
	\begin{equation*}
		\sum_{z\in\mathcal{N}_{\T_{\ell+1}}^*(T)}\norm{\omega_{\text{ref}}^{1/2-s}\lambda_z\left(\varphi_z^{\ell+1}\circ F_T\right)}{L^2(T_{\text{ref}})}^2 \simeq \norm{\omega_{\text{ref}}^{1/2-s}(v_{\ell+1}\circ F_T)}{L^2(T_{\text{ref}})}^2.
	\end{equation*}
    Thus, together with \eqref{proofs:scalingfundamental}, we conclude
	\begin{equation*}
		\sum_{z\in\mathcal{N}_{\T_{\ell+1}}^*}\norm{\widetilde{h}_{\ell}^{-s}\lambda_z\varphi_z^{\ell+1}}{L^2(T)}^2 \simeq h_{\ell}(T)^{d-2s}\norm{\omega_{\text{ref}}^{1/2-s}(v_{\ell+1}\circ  F_T)}{L^2(T_{\text{ref}})}^2
		=\norm{\widetilde{h}^{-s}_{\ell}v_{\ell+1}}{L^2(T)}^2.
	\end{equation*}
	Summing over all $T\in\T_{\ell}$ completes the proof.
\end{fatproof}

Finally, we can present the proof of Theorem~\ref{thm:mainresult1} (b).
\begin{lemma}
	 Let $\T_{\ell}\in\mathbb{T}$, $\mathcal{M}_{\ell}\subseteq \T_{\ell}$ and $\T_{\ell+1}:=\operatorname{refine}(\T_{\ell}, \mathcal{M}_{\ell})$. Let $u_{\ell}$ and $u_{\ell+1}$ be the discrete solutions with respect to the triangulations $\T_{\ell}$ and $\T_{\ell+1}$. Then, there holds
	\begin{align}\label{proofs:discrrelpre}
	\vertiii{u_{\ell+1}-u_{\ell}}^2\lesssim   \sum_{z\in \mathcal{N}_{\T_{\ell+1}}^*\setminus\mathcal{N}_{\T_{\ell}} }\tau_{\ell}\left(\varphi_z^{\ell+1}\right)^2
	\end{align}
	as well as 
	\begin{align}\label{proofs:discrrel}
		\vertiii{u_{\ell+1}-u_{\ell}} \lesssim \tau_{\ell}(\T_{\ell}\setminus\T_{\ell+1}),
	\end{align}
	where the hidden constants depend only on $s$, $d$, $\Omega$ and the $\gamma$-shape regularity of $\T_{\ell}$.
\end{lemma}

\begin{fatproof} Let $I_{\ell}:\S_0^1(\T_{\ell+1})\rightarrow \S_0^1(\T_{\ell})$ be the nodal interpolation operator. As in \cite[Lemma 15]{AFFKP15}, we choose coefficients $\lambda_z$, where $z\in\mathcal{N}_{\T_{\ell+1}}^*\setminus\mathcal{N}_{\T_{\ell}}$, such that
\begin{align*}
	(1-I_{\ell})(u_{\ell+1}-u_{\ell}) = \sum_{z\in \mathcal{N}_{\T_{\ell+1}}^*\setminus\mathcal{N}_{\T_{\ell}}}\lambda_z\varphi_z^{\ell+1}.
\end{align*}
For $z\in \mathcal{N}_{\T_{\ell+1}}^*\setminus\mathcal{N}_{\T_{\ell}}$, let $\mathbb{G}_z:\S_0^1(\T_{\ell+1})\rightarrow\operatorname{span}\{\varphi_z^{\ell+1}\}$ denote the Galerkin-projection onto $\operatorname{span}\{\varphi_z^{\ell+1}\}$.
	With this notation, by using Galerkin-orthogonality and since $\mathbb{G}_z$ is self-adjoint with respect to $a(\cdot, \cdot)$, we obtain
	\begin{equation*}
	\begin{split}
	\vertiii{u_{\ell+1}-u_{\ell}}^2 &= a(u_{\ell+1}-u_{\ell}, (1-I_{\ell})(u_{\ell+1}-u_{\ell}))
	= \sum_{z\in \mathcal{N}_{\T_{\ell+1}}^*\setminus\mathcal{N}_{\T_{\ell}}} a\left(u_{\ell+1}-u_{\ell}, \lambda_z\varphi_z^{\ell+1}\right) \\
	&= \sum_{z\in \mathcal{N}_{\T_{\ell+1}}^*\setminus\mathcal{N}_{\T_{\ell}}} a\left(u_{\ell+1}-u_{\ell}, \lambda_z\mathbb{G}_z\varphi_z^{\ell+1}\right)
	= \sum_{z\in \mathcal{N}_{\T_{\ell+1}}^*\setminus\mathcal{N}_{\T_{\ell}}} a\left(\mathbb{G}_z(u_{\ell+1}-u_{\ell}), \lambda_z\varphi_z^{\ell+1}\right).
	\end{split}
	\end{equation*}
	The Cauchy-Schwarz inequality yields
	\begin{equation}\label{discrreleq}
	\vertiii{u_{\ell+1}-u_{\ell}}^2\leq \Bigg(\sum_{z\in \mathcal{N}_{\T_{\ell+1}}^*\setminus\mathcal{N}_{\T_{\ell}}}\vertiii{\mathbb{G}_z(u_{\ell+1}-u_{\ell})}^2\Bigg)^{\frac{1}{2}}\Bigg(\sum_{z\in \mathcal{N}_{\T_{\ell+1}}^*\setminus\mathcal{N}_{\T_{\ell}}}\vertiii{\lambda_z\varphi_z^{\ell+1}}^2\Bigg)^{\frac{1}{2}}.
	\end{equation}
	
	According to Lemma \ref{normequivalence} there holds
	\begin{align*}
	\vertiii{\varphi_z^{\ell+1}}\lesssim \norm{\widetilde{h}_{\ell}^{-s}\varphi_z^{\ell+1}}{L^2(\Omega)}
	\end{align*}
	for all $z\in \mathcal{N}_{\T_{\ell+1}}^*\setminus\mathcal{N}_{\T_{\ell}}$.
	Together with Lemma \ref{inversetriangle} we conclude
	\begin{align*}
	\sum_{z\in \mathcal{N}_{\T_{\ell+1}}^*\setminus\mathcal{N}_{\T_{\ell}}}\vertiii{\lambda_z\varphi_z^{\ell+1}}^2 &\lesssim \sum_{z\in \mathcal{N}_{\T_{\ell+1}}^*\setminus\mathcal{N}_{\T_{\ell}}}\norm{\widetilde{h}_{\ell}^{-s}\lambda_z\varphi_z^{\ell+1}}{L^2(\Omega)}^2 \lesssim \norm{\sum_{z\in \mathcal{N}_{\T_{\ell+1}}^*\setminus\mathcal{N}_{\T_{\ell}}}\widetilde{h}_{\ell}^{-s}\lambda_z\varphi_z^{\ell+1}}{L^2(\Omega)}^2 \\
	&=\norm{\widetilde{h}_{\ell}^{-s}(1-I_{\ell})(u_{\ell+1}-u_{\ell})}{L^2(\Omega)}^2.
	\end{align*}

	By using Theorem \ref{thm:NodalSZ} and Theorem \ref{scproperties}(b) we get
	\begin{equation*}
	\norm{\widetilde{h}_{\ell}^{-s}(1-I_{\ell})(u_{\ell+1}-u_{\ell})}{L^2(\Omega)}^2\simeq \norm{\widetilde{h}_{\ell}^{-s}(1-J_{\ell})(u_{\ell+1}-u_{\ell})}{L^2(\Omega)}^2 \lesssim \vertiii{u_{\ell+1}-u_{\ell}}^2,
	\end{equation*}
	
	i.e.,
	\begin{equation*}
\sum_{z\in \mathcal{N}_{\T_{\ell+1}}^*\setminus\mathcal{N}_{\T_{\ell}}}\vertiii{\lambda_z\varphi_z^{\ell+1}}^2 \lesssim \vertiii{u_{\ell+1}-u_{\ell}}^2.
	\end{equation*}
	From \eqref{discrreleq} and Remark \ref{setting:Galerkinremark} we deduce
	\begin{equation*}
	\vertiii{u_{\ell+1}-u_{\ell}}^2 \lesssim \sum_{z\in \mathcal{N}_{\T_{\ell+1}}^*\setminus\mathcal{N}_{\T_{\ell}}}\vertiii{\mathbb{G}_z(u_{\ell+1}-u_{\ell})}^2 =   \sum_{z\in \mathcal{N}_{\T_{\ell+1}}^*\setminus\mathcal{N}_{\T_{\ell}}}\tau_{\ell}(\varphi_z^{\ell+1})^2,
	\end{equation*}
	which finishes the proof of \eqref{proofs:discrrelpre}.
	
	In the two-dimensional case it was noticed in \cite[Proof of Theorem 3.2]{PRS20} that, for any node $z\in \mathcal{N}_{\T_{\ell+1}}^*\setminus\mathcal{N}_{\T_{\ell}}$, there holds $\varphi_z^{\ell+1} = \varphi_z^{\ell,u}$, which, in combination with \eqref{proofs:discrrelpre} implies \eqref{proofs:discrrel}.

	The three-dimensional case is more involved: For $z\in \mathcal{N}_{\T_{\ell+1}}^*\setminus\mathcal{N}_{\T_{\ell}}$, we have to compare $\tau_{\ell}\left(\varphi_z^{\ell+1}\right)$ and $\tau_{\ell}\left(\varphi_z^{\ell,u}\right)$. For any node $z\in \mathcal{N}_{\T_{\ell+1}}^*\setminus\mathcal{N}_{\T_{\ell}}$ there exists an element $T\in\T_{\ell}$ such that $z$ is either contained in a face of $T$ or is the midpoint of an edge of $T$. That is, we can write $\mathcal{N}_{\T_{\ell+1}}^*\setminus\mathcal{N}_{\T_{\ell}} = \mathcal{E}_{\T_{\ell+1}}\cup \mathcal{F}_{\T_{\ell+1}}$, where $\mathcal{E}_{\T_{\ell+1}}$ is the set of new nodes that are the midpoint of edges of elements in $\T_{\ell}$ and $\mathcal{F}_{\T_{\ell+1}}$ is the set of new nodes contained in faces of elements in $\T_{\ell}$.
	
	For $z\in \mathcal{F}_{\T_{\ell+1}}$, there holds $\varphi_z^{\ell+1} = \varphi_z^{\ell,u}$ and consequently $\tau_{\ell}\left(\varphi_z^{\ell+1}\right)=\tau_{\ell}\left(\varphi_z^{\ell,u}\right)$, which implies 
	\begin{align}\label{nodalsplit}
		\sum_{z\in \mathcal{F}_{\T_{\ell+1}}}\tau_{\ell}(\varphi_z^{\ell+1})^2\leq \tau_{\ell}\left(\T_{\ell}\setminus\T_{\ell+1}\right)^2.
	\end{align}

	If $z\in \mathcal{E}_{\T_{\ell+1}}$, then there exists an edge $E_z$ of an element $T\in\T_{\ell}$, such that $z$ is the midpoint of $E_z$. 
 	We consider the set of nodes
	\begin{align*}
		\mathcal{N}_{\Gamma_{\ell}(E_z)}:=\left\{z'\in\mathcal{N}^*_{\Tunif}\setminus\mathcal{N}_{\T_{\ell+1}}\ \big|\ z'\ \text{is\ contained\ in\ a\ face\ } F\in\Gamma_{\ell}(E_z)\right\},
	\end{align*}
	where $\Gamma_{\ell}(E_z)$ is the set of all faces $F$ in the triangulation $\T_{\ell}$ with the property that $E_z$ is an edge of $F$.
	We observe that there holds
	\begin{align*}
		\varphi_z^{\ell+1} = \varphi_z^{\ell,u}+\sum_{z'\in\mathcal{N}_{\Gamma_{\ell}(E_z)}} \frac{1}{2}\varphi_{z'}^{\ell,u},
	\end{align*}
	where the factor $1/2$ comes from $\varphi_z^{\ell+1}(z') = 1/2$ for all $z'\in\mathcal{N}_{\Gamma_{\ell}(E_z)}$.
	Galerkin orthogonality shows
	\begin{align*}
		\left|\SCP{f}{\varphi_z^{\ell+1}}{L^2(\Omega)}-a\left(u_{\ell},\varphi_z^{\ell+1}\right)\right|
		&= \left|a\left(u_{\ell,u}-u_{\ell},\varphi_z^{\ell+1}\right)\right| \\
		&\leq \left|a\left(u_{\ell,u}-u_{\ell},\varphi_z^{\ell,u}\right)\right|+ \sum_{z'\in\mathcal{N}_{\Gamma_{\ell}(E_z)}}\left|a\left(u_{\ell,u}-u_{\ell},\varphi_{z'}^{\ell,u}\right)\right|,
	\end{align*}

	 and, for a node $z\in \mathcal{N}_{\T_{\ell+1}}^*\setminus\mathcal{N}_{\T_{\ell}}$, Lemma \ref{normequivalence} (applied to $\T_{\ell+1}$), the pointwise inequality $\varphi_z^{\ell+1}\geq \varphi_z^{\ell,u}$ and Lemma \ref{normequivalence} (applied to $\Tunif$) lead to 
	\begin{align}\label{analog1}
		\vertiii{\varphi_z^{\ell+1}}\gtrsim\norm{\widetilde{h}^{-s}_{\ell}\varphi_z^{\ell+1}}{L^2(\Omega)}\geq \norm{\widetilde{h}^{-s}_{\ell}\varphi_z^{\ell,u}}{L^2(\Omega)}\gtrsim \vertiii{\varphi_z^{\ell,u}}.
	\end{align}
	Furthermore, for $z'\in\mathcal{N}_{\Gamma_{\ell}(E_z)}$, there holds $\varphi_z^{\ell+1}\geq \frac{1}{2}\varphi_{z'}^{\ell,u}$, and with the same arguments as \eqref{analog1} one shows for all $z'\in\mathcal{N}_{\Gamma_{\ell}(E_z)}$
	\begin{align*}
		\vertiii{\varphi_z^{\ell+1}}\gtrsim\vertiii{\varphi_{z'}^{\ell,u}}.
	\end{align*}
	Altogether, for $z\in\mathcal{E}_{\T_{\ell+1}}$, we obtain
	\begin{align*}
		\tau_{\ell}\left(\varphi^{\ell+1}_z\right) &= \frac{\left|\SCP{f}{\varphi^{\ell+1}_z}{L^2(\Omega)}-a\left(u_{\ell},\varphi^{\ell+1}_z\right)\right|}{\vertiii{\varphi^{\ell+1}_z}} \\
		&\lesssim \tau_{\ell}\left(\varphi^{\ell,u}_z\right)+\sum_{z'\in\mathcal{N}_{\Gamma_{\ell}(E_z)}}\tau_{\ell}\left(\varphi_{z'}^{\ell,u}\right). 
	\end{align*}
	
	Since the number of nodes in $\mathcal{N}_{\Gamma_{\ell}(E_z)}$ is bounded by a constant depending only on the $\gamma$-shape regularity of $\T_{\ell}$, we conclude
	\begin{align*}
		\tau_{\ell}\left(\varphi^{\ell+1}_z\right)^2\lesssim \tau_{\ell}\left(\varphi_z^{\ell,u}\right)^2+\sum_{z'\in\mathcal{N}_{\Gamma_{\ell}(E_z)}}\tau_{\ell}\left(\varphi_{z'}^{\ell,u}\right)^2.
	\end{align*}
	Summing over all $z\in\mathcal{E}_{\T_{\ell+1}}$ yields
	
	\begin{align*}
		\sum_{z\in\mathcal{E}_{\T_{\ell+1}}}\tau_{\ell}\left(\varphi_z^{\ell+1}\right)^2 &\lesssim \tau_{\ell}\left(\T_{\ell}\setminus\T_{\ell+1}\right)^2+\sum_{z\in\mathcal{E}_{\T_{\ell+1}}}\sum_{z'\in\mathcal{N}_{\Gamma_{\ell}(E_z)}}\tau_{\ell}\left(\varphi_{z'}^{\ell,u}\right)^2 \\
		&\lesssim \tau_{\ell}\left(\T_{\ell}\setminus\T_{\ell+1}\right)^2,
	\end{align*}
	where the last inequality is due to the fact that every $z'\in\mathcal{N}_{\Gamma_{\ell}(E_z)}$ is counted at most three times and that every element $T\in\T_{\ell}$, which has a face $F$ belonging to $\Gamma_{\ell}(E_z)$, is necessarily in $\T_{\ell}\setminus\T_{\ell+1}$. Thus, together with \eqref{nodalsplit} we have 
	\begin{align*}
		\sum_{z\in \mathcal{N}_{\T_{\ell+1}}^*\setminus\mathcal{N}_{\T_{\ell}}}\tau_{\ell}(\varphi_z^{\ell+1})^2 &= \sum_{z\in \mathcal{F}_{\T_{\ell+1}}}\tau_{\ell}(\varphi_z^{\ell+1})^2+\sum_{z\in \mathcal{E}_{\T_{\ell+1}}}\tau_{\ell}(\varphi_z^{\ell+1})^2 \\
		&\lesssim \tau_{\ell}\left(\T_{\ell}\setminus\T_{\ell+1}\right)^2,
	\end{align*}
	
	which, together with \eqref{proofs:discrrelpre}, proves \eqref{proofs:discrrel} in 3D.
\end{fatproof}

\subsection{Proof of Theorem \ref{thm:mainresult1} (c)}

The proof of Theorem \ref{thm:mainresult1} (c) requires a modification of Lemma \ref{proofs:hequivalence}.

\begin{lemma}\label{modhequivalence}
	 Assume $d=2$. Let $\T_{L}\in\mathbb{T}$ with uniform refinement $\T_{L,u} := \operatorname{refine}(\T_{L},\T_{L})$, as well as a triangulation $\T_{\ell}\in\operatorname{refine}(\T_{L})$ with uniform refinement $\T_{\ell,u} := \operatorname{refine}(\T_{\ell},\T_{\ell})$ and consider an element $T\in\T_{L}\cap\T_{\ell}$. Then, there holds $\mathcal{N}_{\T_{L,u}}^*(T)\setminus\mathcal{N}_{\T_L} = \mathcal{N}_{\T_{\ell,u}}^*(T)\setminus\mathcal{N}_{\T_{\ell}}$, and, for any node $z\in\mathcal{N}_{\T_{L,u}}^*(T)\setminus\mathcal{N}_{\T_{L}}$, the associated nodal basis function $\varphi_z^{L,u}$ is the same in $\S_0^1(\T_{L,u})$ and $\S_0^1(\T_{\ell,u})$. Furthermore, there holds
	 
	\begin{equation}\label{proofs:Llequival}
		\norm{\widetilde{h}_{\ell}^{-s}\varphi_z^{L,u}}{L^2(\Omega)}\simeq \norm{\widetilde{h}_{L}^{-s}\varphi_z^{L,u}}{L^2(\Omega)},
	\end{equation}
	where the hidden constants depend only on $s$, $d$, and the $\gamma$-shape regularity of $\T_{L}$.
\end{lemma}

\begin{fatproof}
	In 2D, it is easy to check that the refinement of a marked element does not depend on its neighbors, and in combination with $T\in\T_L\cap\T_{\ell}$, this implies $\mathcal{N}_{\T_{L,u}}^*(T)\setminus\mathcal{N}_{\T_L} = \mathcal{N}_{\T_{\ell,u}}^*(T)\setminus\mathcal{N}_{\T_{\ell}}$. Furthermore, in 2D the nodal basis function $\varphi_z^{L,u}$ associated to a node $z\in\mathcal{N}_{\T_{L,u}}^*(T)\setminus\mathcal{N}_{\T_{L}}$ depends only on $T$. Therefore $\varphi_z^{L,u}$ is the same in $\S_0^1(\T_{L,u})$ and $\S_0^1(\T_{\ell,u})$.
	
    It remains to show \eqref{proofs:Llequival}. Due to the pointwise inequality $\widetilde{h}_{\ell}^{-s}\geq \widetilde{h}_L^{-s}$, the upper inequality is clear and we only have to show the lower estimate.
    
	Consider a node $z\in\mathcal{N}_{\T_{L,u}}^*(T)\setminus\mathcal{N}_{\T_{L}}$. Then, $z$ is the midpoint of an edge $E_z$ of $T$. Let $T'\in\T_{\ell}$ be an element in the fine mesh $\T_{\ell}$, which shares $E_z$. Then, $T'$ is either itself an element of the coarse mesh $\T_L$ or it is a son of an element $T^*\in\T_L$. Applying one or two NVB bisection steps (depending on $T'\in\T_L$ or $T'$ being a son of an element in $\T_L$) then leads to elements in the uniform refinement of the coarse mesh $\T_{L,u}$.
	This means that on $\operatorname{supp}(\varphi_z^{L,u})$ the uniform refinement is locally finer than $\T_{\ell}$.

	In terms of the corresponding weight-functions $\widetilde{h}_{\ell}$ and $\widetilde{h}_{L,u}$ (where $\widetilde{h}_{L,u}$ is the weight-function associated to $\T_{L,u}$), this implies that $\widetilde{h}_{L,u}^ {-s}\geq \widetilde{h}_{\ell}^ {-s}$ pointwise on $\operatorname{supp}(\varphi_z^{L,u})$. From Lemma \ref{proofs:hequivalence} we infer 
	\begin{align*}
		\left\|\widetilde{h}_{L,u}^{-s}\varphi_z^{L,u}\right\|_{L^ 2(\Omega)}\simeq \left\|\widetilde{h}_{L}^{-s}\varphi_z^{L,u}\right\|_{L^ 2(\Omega)},
	\end{align*}
	which proves the lower estimate in \eqref{proofs:Llequival}.
\end{fatproof}

Finally, we prove stability of the two-level estimator in two dimensions.

\begin{lemma}\label{modtheorem}
	 Assume $d=2$. Let $\T_{L}\in\mathbb{T}$ with uniform refinement $\T_{L,u} := \operatorname{refine}(\T_{L},\T_{L})$, as well as $\T_{\ell}\in\operatorname{refine}(\T_{L})$. Furthermore, let $u_{L}$ and $u_{\ell}$ be the discrete solutions associated to $\T_{L}$ and $\T_{\ell}$. Then,
	\begin{equation}\label{prediscrstab}
	\sum_{T\in\T_{L}\cap \T_{\ell}}\ \sum_{z\in\mathcal{N}_{\T_{L,u}}^*(T)\setminus\mathcal{N}_{\T_{L}}}\frac{\left|a\left(u_{L}-u_{\ell}, \varphi_z^{L,u}\right)\right|^2}{\vertiii{\varphi_z^{L,u}}^2} \lesssim \vertiii{u_{L}-u_{\ell}}^2,
	\end{equation}
    as well as 
    \begin{align}\label{discrstab}
    	\big|\tau_{\ell}(\T_L\cap\T_{\ell})-\tau_L(\T_L\cap\T_{\ell})\big|\lesssim \vertiii{u_L-u_{\ell}},
    \end{align}
	where the hidden constants depend only on $s$, $d$, and the $\gamma$-shape regularity of $\T_{L}$.
	
\end{lemma}

\begin{fatproof}
	Let $T\in\T_{L}\cap\T_{\ell}$ be arbitrary and $z\in\mathcal{N}_{\T_{L,u}}^*(T)\setminus\mathcal{N}_{\T_{L}}$. Due to $\operatorname{supp}(\varphi_z^{L,u})\subseteq \Omega_L[T]$ and Lemma \ref{modhequivalence}, analogous arguments as in the proof of Lemma \ref{proofs:localestimate} show
	\begin{equation*}
	\left|a\left(u_{L}-u_{\ell},\varphi_z^{L,u}\right)\right|
	\lesssim \norm{\widetilde{h}_{\ell}^s(-\Delta)^s(u_{L}-u_{\ell})}{L^2(\Omega_{L}[T])}\vertiii{\varphi_z^{L,u}}.
	\end{equation*}
	
	Consequently, summing over all $z\in \mathcal{N}_{\T_{L,u}}^*(T)\setminus\mathcal{N}_{\T_{L}}$, and since $\#\{\mathcal{N}_{\T_{L,u}}^*(T)\setminus\mathcal{N}_{\T_{L}}\}$ is bounded with an upper bound depending only on $d$, we conclude
	\begin{equation*}
	\begin{split}
	\sum_{z\in \mathcal{N}_{\T_{L,u}}^*(T)\setminus\mathcal{N}_{\T_{L}}}\frac{\left|a\left(u_{L}-u_{\ell},\varphi_z^{L,u}\right)\right|^2}{\vertiii{\varphi_z^{L,u}}^2}
	&\lesssim \sum_{z\in \mathcal{N}_{\T_{L,u}}^*(T)\setminus\mathcal{N}_{\T_{L}} }\norm{\widetilde{h}_{\ell}^s(-\Delta)^s(u_{L}-u_{\ell})}{L^2(\Omega_{L}[T])}^2 \\
	&\lesssim \norm{\widetilde{h}_{\ell}^s(-\Delta)^s(u_{L}-u_{\ell})}{L^2(\Omega_{L}[T])}^2.
	\end{split}
	\end{equation*}
	Summing over $T\in\T_{L}\cap\T_{\ell}$, using $\gamma$-shape regularity and Theorem \ref{inversestimate}(b) then shows 
	\begin{equation*}
	\begin{split}
	\sum_{T\in\T_{L}\cap \T_{\ell}}\ \sum_{z\in\mathcal{N}_{\T_{L,u}}^*(T)\setminus\mathcal{N}_{\T_{L}}}\frac{\left|a\left(u_{L}-u_{\ell}, \varphi_z^{L,u}\right)\right|^2}{\vertiii{\varphi_z^{L,u}}^2} &\lesssim \sum_{T\in\T_{L}\cap \T_{\ell}}\norm{\widetilde{h}_{\ell}^s(-\Delta)^s(u_{L}-u_{\ell})}{L^2(\Omega_{L}[T])}^2 \\ 
	&\lesssim \norm{\widetilde{h}_{\ell}^s(-\Delta)^s(u_{L}-u_{\ell})}{L^2(\Omega)}^2 \\
	&\lesssim \vertiii{u_L-u_{\ell}}^2,
	\end{split}
	\end{equation*}
	which is \eqref{prediscrstab}. As in \cite[Proof of Theorem 3.2]{PRS20}, the reverse triangle inequality and the equality $\varphi_z^{L,u}=\varphi_z^{\ell,u}$ (see Lemma \ref{modhequivalence}) show
	\begin{align*}
			\big|\tau_{L}(\T_L\cap\T_{\ell})-\tau_{\ell}(\T_L\cap\T_{\ell})\big| &= \Bigg|\left(\sum_{T\in\T_{L}\cap\T_{\ell}}\tau_{L}(T)^2\right)^{1/2}-\left(\sum_{T\in\T_{L}\cap\T_{\ell}}\tau_{\ell}(T)^2\right)^{1/2}\Bigg| \\ 
			&\leq \left(\sum_{T\in\T_L\cap\T_{\ell}}\ \sum_{z\in\mathcal{N}_{\widehat{\T}_{L,u}}^*(T)\setminus\mathcal{N}_{\T_{L}}}\frac{\left|a\left(u_{L}-u_{\ell}, \varphi_z^{L,u}\right)\right|^2}{\vertiii{\varphi_z^{L,u}}^2}\right)^{1/2} \\
			&\lesssim \vertiii{u_L-u_{\ell}},
	\end{align*}
	which finishes the proof.
	
\end{fatproof}

\begin{remark}
	For the proof of Lemma \ref{modtheorem} it is essential that, for any node $z\in\mathcal{N}_{\T_{L,u}}^*(T)\setminus\mathcal{N}_{\T_{L}}$, the associated nodal basis function $\varphi_z^{L,u}$ is the same in $\S_0^1(\T_{L,u})$ and $\S_0^1(\T_{\ell,u})$. In the case $d=3$ this fact is not longer true. That is, the proof of Lemma \ref{modtheorem} cannot directly be generalized to $d=3$ and stability (and consequently also convergence with optimal algebraic rates) of the two-level error estimator in three dimensions still remains open.
\end{remark}

	\section{Numerical experiments}

We illustrate our work by carrying out some numerical experiments. We consider the model problem in its weak form \eqref{setting:weakform}
in dimension $d=2$. As in \cite{FMP21a}, we choose $\Omega$ to be either the unit circle or a L-shaped domain.

\subsection{Aspects of our implementation}

We implemented Algorithm \ref{SEMR} in MATLAB R2022a. Regarding the steps SOLVE, ESTIMATE and REFINE (steps 1., 2., 4. in Algorithm \ref{SEMR}), we note:
\begin{itemize}
	\item For a given triangulation $\T_{\ell}\in\mathbb{T}$, the computation of the associated discrete solution $u_{\ell}\in\S_0^1(\T_{\ell})$ was done by using the already existing MATLAB code from \cite{FEM1}. The unbounded domain $\R^2$ is replaced by a circle around the domain $\Omega$ and the integrals occurring in the finite element discretization are transformed by a Duffy transformation and then computed by quadrature formulas.
	
	\item Computing the two-level error estimator \eqref{setting:estimator} requires the evaluation of the expressions $a(u_{\ell},\varphi_z^{\ell,u})$ and $a(\varphi_z^{\ell,u}, \varphi_z^{\ell,u})$. Again, these evaluations were computed by using the MATLAB code from \cite{FEM1} and quadrature formulas were used to compute $\SCP{f}{\varphi_z^{\ell,u}}{L^2(\Omega)}$.
	
	\item Based on a set of marked elements $\mathcal{M}_{\ell}\subseteq \T_{\ell}$, the refined mesh $\T_{\ell+1} := \operatorname{refine}(\T_{\ell},\mathcal{M}_{\ell})$ is obtained by NVB. Recall that we assumed that every marked element $T\in\mathcal{M}_{\ell}$ is bisected three times into four sons. In order to compute refined meshes we used an existing MATLAB code from \cite{FEM2}, which generates NVB refinements with the property that every marked element is split into four sons.
\end{itemize}

\subsection{Unit circle with constant right-hand side}

For the first example we choose $\Omega$ to be the unit circle, i.e., $\Omega := B_1(0)$ and set $f := 2^ {2s}\Gamma(1+s)^2$. With this domain and right-hand side, the exact solution to \eqref{setting:weakform} is known (see, e.g., \cite{exactsol}) and is given by $u(x) := (1-|x|^2)^s_+$, where $g_+:=\max\{g,0\}$. The energy norm of the exact solution $u$ can easily be computed by 
\begin{align*}
	a(u,u) = \int_{B_1(0)}fu\ \diff x = 2^{2s}\Gamma(1+s)^2\frac{2\pi}{2s+2}.
\end{align*}

\begin{figure}[h]
	\begin{minipage}{0.5\textwidth}
		\begin{tikzpicture}
  \begin{loglogaxis}[
    width=1\textwidth, height=6cm,     
    grid = major,
    grid style={dashed, gray!30},
    axis background/.style={fill=white},
    xlabel=DOFs N,
    legend style={at={(0.02,0.02)},anchor=south west,font=\footnotesize},
    xtick = {},
    ]
      
    \addplot+[solid,mark=square,mark size=2pt,mark options={line width=1.0pt}] table    
    [
    x index = 0,
    y index = 1,
    ]{ex_25_30_8300_circle_adaptive.txt};    
    \addlegendentry{err adap}

    \addplot+[solid,mark=star,mark size=2pt,mark options={line width=1.0pt}] table    
    [
    x index=0,
    y index = 2,
    ]{ex_25_30_8300_circle_adaptive.txt};    
    \addlegendentry{est adap}

    \addplot+[solid,mark=triangle,mark size=2pt,mark options={line width=1.0pt}] table    
    [
    x index=0,
    y index=1,
    ]{ex_25_30_8300_circle_uniform.txt};    
    \addlegendentry{err unif}
    
    \addplot+[solid,mark=diamond,mark size=2pt,mark options={line width=1.0pt}] table    
    [
    x index=0,
    y index=2,
    ]{ex_25_30_8300_circle_uniform.txt};    
    \addlegendentry{est unif}

    \addplot [black,dashed,mark options={line width=1.0pt}] expression [domain=120:8000, samples=15] {7*x^(-1/2)} node [below, xshift=-1cm, yshift=0.3cm] {$N^{-1/2}$};

    \addplot [black,dashed ] expression [domain=120:8000, samples=15] {2.9*x^(-1/4)} node [above,yshift=0.25cm,xshift=-0.1cm] {$N^{-1/4}$};
  \end{loglogaxis} 
\end{tikzpicture}
	\end{minipage}
	\begin{minipage}{0.5\textwidth}
		\begin{tikzpicture}
  \begin{loglogaxis}[
    width=1\textwidth, height=6cm,     
    grid = major,
    grid style={dashed, gray!30},
    axis background/.style={fill=white},
    xlabel=DOFs N,
    legend style={at={(0.02,0.02)},anchor=south west,font=\footnotesize},
    xtick = {},
    ]
      
    \addplot+[solid,mark=square,mark size=2pt,mark options={line width=1.0pt}] table    
    [
    x index = 0,
    y index = 1,
    ]{ex_75_30_7000_circle_adaptive_constant.txt};    
    \addlegendentry{err adap}

    \addplot+[solid,mark=star,mark size=2pt,mark options={line width=1.0pt}] table    
    [
    x index=0,
    y index = 2,
    ]{ex_75_30_7000_circle_adaptive_constant.txt};    
    \addlegendentry{est adap}

    \addplot+[solid,mark=triangle,mark size=2pt,mark options={line width=1.0pt}] table    
    [
    x index=0,
    y index=1,
    ]{ex_75_60_8300_circle_uniform.txt};    
    \addlegendentry{err unif}
    
    \addplot+[solid,mark=diamond,mark size=2pt,mark options={line width=1.0pt}] table    
    [
    x index=0,
    y index=2,
    ]{ex_75_60_8300_circle_uniform.txt};    
    \addlegendentry{est unif}

    \addplot [black,dashed,mark options={line width=1.0pt}] expression [domain=120:8000, samples=15] {7*x^(-1/2)} node [below, xshift=-1cm, yshift=0.3cm] {$N^{-1/2}$};

    \addplot [black,dashed ] expression [domain=120:8000, samples=15] {2.3*x^(-1/4)} node [above,yshift=0.25cm,xshift=-0.1cm] {$N^{-1/4}$};
  \end{loglogaxis} 
\end{tikzpicture}
	\end{minipage}	
\caption{Error and two-level estimator for a constant right-hand side on the unit circle for uniform and adaptive mesh refinement. Left: $s=0.25$ and $\theta = 0.3$. Right: $s=0.75$ and $\theta = 0.3$.}
\label{fig:ex11}
\end{figure}
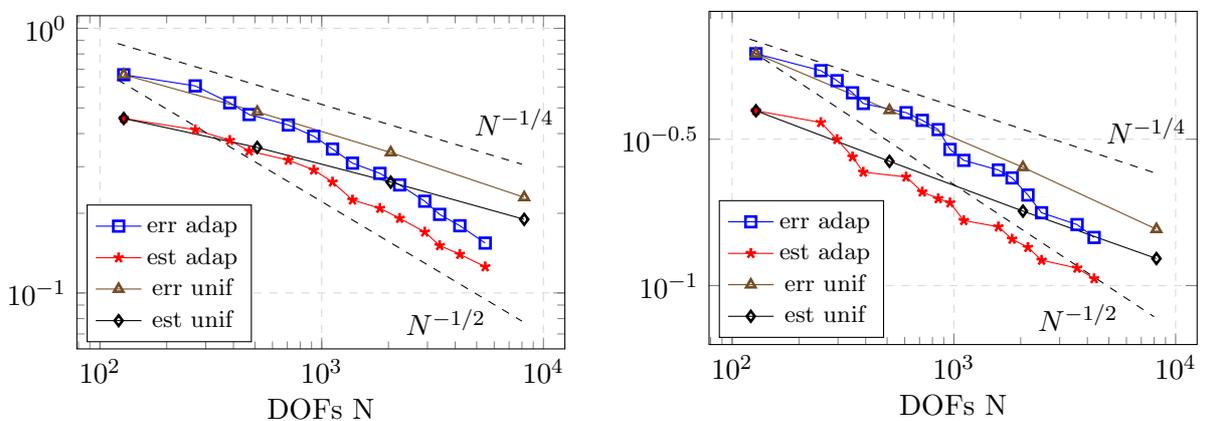

Figure \ref{fig:ex11} shows the uniform and adaptive error estimators, as well as the uniform and adaptive errors for the cases $s=0.25$ and $\theta=0.3$, as well as $s=0.75$ and $\theta=0.3$.

We observe that uniform refinement leads to a reduced rate of $N^{-1/4}$ for error and error estimator. The adaptive algorithm, however, leads to the optimal rate of $N^{-1/2}$.

\subsection{L-shaped domain with constant right-hand side}

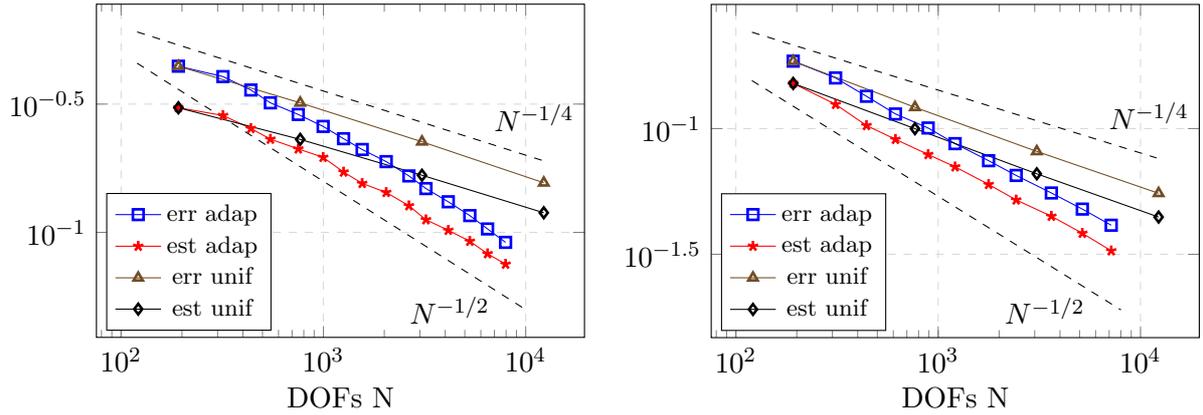
\begin{figure}[h]
	\begin{minipage}{0.5\textwidth}
		\begin{tikzpicture}
  \begin{loglogaxis}[
    width=1\textwidth, height=6cm,     
    grid = major,
    grid style={dashed, gray!30},
    axis background/.style={fill=white},
    xlabel=DOFs N,
    legend style={at={(0.02,0.02)},anchor=south west,font=\footnotesize},
    xtick = {},
    ]
      
    \addplot+[solid,mark=square,mark size=2pt,mark options={line width=1.0pt}] table    
    [
    x index = 0,
    y index = 1,
    ]{ex_25_30_8300_Lshape_adaptive.txt};    
    \addlegendentry{err adap}

    \addplot+[solid,mark=star,mark size=2pt,mark options={line width=1.0pt}] table    
    [
    x index=0,
    y index = 2,
    ]{ex_25_30_8300_Lshape_adaptive.txt};    
    \addlegendentry{est adap}

    \addplot+[solid,mark=triangle,mark size=2pt,mark options={line width=1.0pt}] table    
    [
    x index=0,
    y index=1,
    ]{ex_25_30_12500_Lshape_uniform.txt};    
    \addlegendentry{err unif}
    
    \addplot+[solid,mark=diamond,mark size=2pt,mark options={line width=1.0pt}] table    
    [
    x index=0,
    y index=2,
    ]{ex_25_30_12500_Lshape_uniform.txt};    
    \addlegendentry{est unif}

    \addplot [black,dashed,mark options={line width=1.0pt}] expression [domain=120:10000, samples=15] {5*x^(-1/2)} node [below, xshift=-1cm, yshift=0.3cm] {$N^{-1/2}$};

    \addplot [black,dashed ] expression [domain=120:12000, samples=15] {2*x^(-1/4)} node [above,yshift=0.25cm,xshift=-0.1cm] {$N^{-1/4}$};
  \end{loglogaxis} 
\end{tikzpicture}
	\end{minipage}
	\begin{minipage}{0.5\textwidth}
		\begin{tikzpicture}
  \begin{loglogaxis}[
    width=1\textwidth, height=6cm,     
    grid = major,
    grid style={dashed, gray!30},
    axis background/.style={fill=white},
    xlabel=DOFs N,
    legend style={at={(0.02,0.02)},anchor=south west,font=\footnotesize},
    xtick = {},
    ]
      
    \addplot+[solid,mark=square,mark size=2pt,mark options={line width=1.0pt}] table    
    [
    x index = 0,
    y index = 1,
    ]{ex_75_40_8500_Lshape_adaptive.txt};    
    \addlegendentry{err adap}

    \addplot+[solid,mark=star,mark size=2pt,mark options={line width=1.0pt}] table    
    [
    x index=0,
    y index = 2,
    ]{ex_75_40_8500_Lshape_adaptive.txt};    
    \addlegendentry{est adap}

    \addplot+[solid,mark=triangle,mark size=2pt,mark options={line width=1.0pt}] table    
    [
    x index=0,
    y index=1,
    ]{ex_75_60_12500_Lshape_uniform.txt};    
    \addlegendentry{err unif}
    
    \addplot+[solid,mark=diamond,mark size=2pt,mark options={line width=1.0pt}] table    
    [
    x index=0,
    y index=2,
    ]{ex_75_60_12500_Lshape_uniform.txt};    
    \addlegendentry{est unif}

    \addplot [black,dashed,mark options={line width=1.0pt}] expression [domain=120:8000, samples=15] {1.7*x^(-1/2)} node [below, xshift=-1cm, yshift=0.3cm] {$N^{-1/2}$};

    \addplot [black,dashed ] expression [domain=120:12000, samples=15] {0.8*x^(-1/4)} node [above,yshift=0.25cm,xshift=-0.1cm] {$N^{-1/4}$};
  \end{loglogaxis} 
\end{tikzpicture}
	\end{minipage}	
	\caption{Error and two-level estimator for a constant right-hand side on the L-shaped domain for uniform and adaptive mesh refinement. Left: $s=0.25$ and $\theta = 0.3$. Right: $s=0.75$ and $\theta = 0.4$.}
	\label{fig:ex13}
\end{figure}

In our second numerical example we choose $\Omega := (-1,1)^2\setminus[0,1)^2$ and $f\equiv 1$. The exact solution of this problem is unknown, therefore we have to extrapolate the energy of the exact solution from the energies of the discrete solutions. Figure \ref{fig:ex13} shows the computed errors and estimators for $s=0.25$ and $\theta=0.3$, as well as $s=0.75$ and $\theta = 0.4$, respectively. 

Qualitatively, we observe the same as in the first example: Uniform refinement yields the reduced rate $N^{-1/4}$, whereas adaptive refinement leads to improved rates of convergence.
\vspace{1cm}

	\bibliographystyle{amsalpha}
	\bibliography{bibliog}
\end{document}